\documentclass{article}
\usepackage{graphicx}
\usepackage{amsmath}
\usepackage{amsfonts}
\usepackage{epsfig}
\newtheorem{proposition}{Proposition}

\renewcommand{\ldots}{\dotsc}

\title{Nonclassical equivalence transformations associated with \\a parameter identification problem
\thanks{This work was partially supported by the Austrian Science
Foundation FWF, Project SFB 1308 "Large Scale Inverse Problems".}}

\author{Nicoleta B\^{i}l\u{a} \thanks{Department of Mathematics and Computer
Science, Fayetteville State University, Fayetteville, NC 28301, USA.
(nbila@uncfsu.edu).} \and Jitse Niesen \thanks{Department of Applied
Mathematics, School of Mathematics, University of Leeds, Leeds LS2
9JT, United Kingdom. (jitse@maths.leeds.ac.uk).}}

\date{}
\begin{document}

\maketitle

\begin{abstract}
A special class of symmetry reductions called nonclassical
equivalence transformations is discussed in connection to a class of
parameter identification problems represented by partial
differential equations. These symmetry reductions relate the forward
and inverse problems, reduce the dimension of the equation, yield
special types of solutions, and may be incorporated into the
boundary conditions as well. As an example, we discuss the nonlinear
stationary heat conduction equation and show that this approach
permits the study of the model on new types of domains. Our
\texttt{MAPLE} routine \texttt{GENDEFNC} which uses the package
\texttt{DESOLV} (authors Carminati and Vu) has been updated for this
propose and its output is the nonlinear partial differential
equation system of the determining equations of the nonclassical
equivalence transformations.
\end{abstract}

\medskip
\par
\noindent AMS: 58J70, 70G65, 35R30, 35G30


\section{Introduction}
One of the fastest developing research fields in the last few years
is the area of inverse problems. These problems arise frequently in
engineering, mathematics, and physics. In particular, parameter
identification problems deal with the identification of physical
parameters from observations of the evolution of a system and
especially arise when the physical laws governing the processes are
known, but the information about the parameters occurring in
equations is needed. In general, these are ill-posed problems, in
the sense that they do not fulfill Hadamard's postulates for all
admissible data: a solution exists, the solution is unique, and the
solution depends continuously on the given data. Arbitrary small
changes in data may lead to arbitrary large changes in the solution
\cite{engl1}.

Consider the following class of parameter identification problems
modeled by partial differential equations (PDEs) of the form
\begin{equation}\label{ip}
F(x, w^{(m)}, E^{(n)})=0,
\end{equation}
where the unknown function $E=E(x)$ is called \textit{parameter},
and, respectively, the arbitrary function $w=w(x)$ is called
\textit{data}, with $x=(x_1,\ldots, x_p)\in \Omega \subset R^p$ a
given domain, where $w^{(m)}$ denotes the function $w$ together with
its partial derivatives up to order~$m$. Suppose that the parameters
and data are analytical functions. The PDE (\ref{ip}) sometimes
augmented with certain boundary conditions is called \textit{the
inverse problem} associated with a \textit{direct (or
\textit{forward}) problem}. The direct problem is the same equation
but the unknown function is the data (for which certain boundary
conditions are imposed as well).

On the other hand, symmetry analysis theory has been widely used to
study nonlinear PDEs. A remarkable number of mathematical physics
models have been successfully analyzed from this point of view.
There is a considerable body of literature on this topic (see, for
example, \cite{Ames}, \cite{bluman1}--\cite{BruGand},
\cite{clarkson1}, \cite{clarkson2}, \cite{GGL1}--\cite{ibragimov},
\cite{meleshko1}--\cite{ovsyannikov1},
\cite{handbookPZ2}--\cite{torrisi3} and references from there).
Sophus Lie ~\cite{lie1} introduced the notion of continuous groups
of transformations (known today as Lie groups of transformations)
and, subsequently, a method for finding the symmetry group
associated with a PDE. This (local) group of transformations acting
on the space of the independent variables and the space of the
dependent variables of the equation has the property that it leaves
the set of all analytical solutions of the PDE invariant. Moreover,
the form or the equation remains unchanged and group-invariant
solutions can be found. Lie's method has been applied extensively to
various mathematical models represented by PDEs. Consequently, new
methods for seeking explicit solutions that cannot be obtained by
Lie's method have been proposed over the years. The nonclassical
method introduced by Bluman and Cole~\cite{bluman1} is a technique
for determining the \textit{nonclassical symmetries} (or
\textit{conditional symmetries}) related to a PDE. These are
transformations that act on the space of the independent variables
and the space of the dependent variables with the property that they
leave only a subset of the set of all analytical solutions
invariant. Knowledge of these classical and nonclassical symmetries
allows one to reduce the dimension of the studied PDE and to
determine solutions of special type. It may happen that certain
group-invariant solutions cannot be found explicitly, but at least
additional information about the studied model can be obtained in
this case. Any classical symmetry is a nonclassical symmetry but not
conversely. Ovsyannikov ~\cite{ovsyannikov1} introduced the notion
of equivalence transformations for a PDE depending on an arbitrary
function and these transformations have been recently generalized by
Meleshko in ~\cite{meleshko1}. The (generalized) equivalence
transformations are groups of transformations acting on the space of
the independent variables, the space of the dependent variables, and
the space of the arbitrary functions with the property that they
leave the differential structure of the equation unchanged. By
dropping the auxiliary conditions that characterize the functional
dependence, Torrisi and Tracin\`{a} introduced the notion of weak
equivalence transformations ~\cite{torrisi3}. These transformations
have been later generalized by Romano and Torrisi in
\cite{torrisi1}. Notice that all these methods do not take into
account the boundary conditions attached to a PDE. Moreover, in the
equivalence transformations approach the arbitrary functions may
depend on the dependent variable only. In our case, for parameter
identification problems of the form (\ref{ip}), the data does not
depend on the parameter -- both data and parameter act on the space
of the independent variables.

The aim of this paper is to continue the study of the parameter
identification problems (\ref{ip}) from the point of view of
symmetry analysis theory. This is a second paper in which we analyze
symmetry reductions suitable for (\ref{ip}). Investigation of
special Lie groups of transformations related to a model may allow
us to reduce its dimension and further obtain information about its
solutions. As discussed in~\cite{bilacar}, finding the classical Lie
symmetries related to the inverse problem (\ref{ip}) might be a very
difficult task and the success of the method depends on the
nonlinearity of the equation. Therefore, considering both the
parameter and data as unknown functions in (\ref{ip}), the classical
Lie method can be applied with less difficulties (in ~\cite{bilacar}
a mathematical model arising in car windshield design was analyzed
using this assumption). Throughout this paper, similar to
\cite{bilacar}, we shall refer to these transformations as
equivalence transformations, taking into account that $w$ was
initially an arbitrary function (these transformations are different
than the equivalence transformations introduced by Ovsyannikov, but
in the literature they are also referred to as equivalence
transformations). Furthermore, once the equivalence transformations
related to (\ref{ip}) have been discussed, it is natural to
investigate the nonclassical symmetries associated with this
equation as well. These new symmetry reductions also relate the
forward and inverse problems as the equivalence transformations do.
To the best of our knowledge, these type of transformations have not
been studied so far for parameter identification problems. Since the
differential structure of the equation is preserved and $w$ is an
arbitrary function, we shall call these transformations
\textit{nonclassical equivalence transformations}. While the
equivalence transformations are found by integrating a linear PDE
system, the nonclassical equivalence transformations are solutions
of a nonlinear PDE system. Typically, a large amount of calculations
are required and, hence, in both cases, the use of a symbolic
manipulation program would represent a great advantage. Therefore,
for this purpose, we have updated our Maple package
\texttt{GENDEFNC} ~\cite{bilaniesen1} (which uses the package
\texttt{DESOLV} by Carminati and Vu ~\cite{desolv}). The
\texttt{GENDEFNC} output is the nonlinear partial differential
equation system of the determining equations of the nonclassical
equivalence transformations. Notice that \texttt{GENDEFNC} is based
on a new algorithm for finding the nonclassical symmetries (see
~\cite{bilaniesen1} for further details). Recently, Bruz\'{o}n and
Gandarias have extended this algorithm to a different case (see
~\cite{BruGand} for more information).

As an example, we shall consider a mathematical model arising in
heat conduction, namely, the nonlinear stationary heat conduction
equation given by
\begin{equation}
\label{eq1} - \hbox{div} (E(x,y) \nabla w(x,y) )=1 \quad
\hbox{in}\quad \Omega,
\end{equation}
where the unknown function $E=E(x,y)$ is the \textit{parameter}, the
arbitrary function  $w=w(x,y)$ is the \textit{data}, and $(x,y)\in
\Omega $ with $\Omega \subset R^2 $ a bounded domain (here $\nabla w
=(w_x,w_y)$ denotes the gradient of $w$). The data function must
also satisfy the Dirichlet boundary condition
\begin{equation}
\label{cond1}
w\vert _{\partial \Omega }=0.
\end{equation}
In 3D, the above problem is related to the heat conduction in a
material occupying a domain $\Omega $ whose temperature is kept zero
at the boundary ~\cite{engl1}. After sufficiently long time, the
temperature distribution $w$  can be modeled by
\begin{equation}
\label{eq3D} - \hbox{div} (E (x,y,z)\nabla w(x,y,z))=f(x,y,z) \quad
\hbox{in} \quad \Omega ,
\end{equation}
where $E$ is the heat conductivity and $f$ represents the heat
sources. For given $E$ and $f$, the forward problem is to find the
temperature distribution $w$ satisfying (\ref{eq3D}) and
(\ref{cond1}). Conversely, the inverse problem is to determine $E$
from (\ref{eq3D}) and (\ref{cond1}) when $w$ is known. While the
direct problem is an elliptic PDE for $w$, the inverse problem is a
linear PDE (with variable coefficients) for $E$. On the other hand,
in the inverse problems approach, if the solution of the forward
problem is unique for each parameter $E$, the parameter-to-solution
map associates with each parameter $E$ the forward problem solution.
Since for the above problem the parameter-to-output map is
nonlinear, (\ref{eq3D}) and (\ref{cond1}) is a nonlinear problem
from this point of view. In addition, when $w=w(E)$, a new dependent
variable can be introduced by using the Kirchoff transformation
~\cite[p.~\mbox{113}]{zwillinger}; if the heat conductivity $E=1$,
then the equation (\ref{eq3D}) becomes Poisson's equation
\cite[p.~\mbox{316}]{Eriksson1}. For simplicity, in this paper, we
discuss the 2D case with the heat sources $f=1$. Since (\ref{eq1})
can be written as
\begin{equation}
\label{eq2} w_xE_x+w_yE_y+E(w_{xx}+w_{yy})=-1,
\end{equation}
it can be easily seen that the vector field
\[
X=w_x\partial _x+w_y\partial
_y-\left[1+E(w_{xx}+w_{yy})\right]\partial _E,
\]
generates the classical Lie symmetries associated with (\ref{eq2})
considered as a PDE with $w$ given and $E$ unknown.

The paper is organized as follows. The equivalence transformations
related to (\ref{eq2}) are presented in \S 2 and the nonclassical
equivalence transformations are discussed in \S 3. We show that the
nonclassical symmetries are related to the Monge equation
(\ref{eqnMonge}), the Monge-Amp\`{e}re equation (\ref{MAinV}), and
the Abel ordinary differential equations (ODEs) of second kind
(\ref{Abel10}) and (\ref{Abel2}). The new symmetry reductions
related to (\ref{eq2}) are given by (\ref{case221}),
(\ref{case222}), (\ref{prop412a}), (\ref{eqncase30}), and,
respectively, (\ref{case3}) and (\ref{eqncase4}) excepting the case
when $A(x,y)=(k_1-k_3y+k_4x)/(k_2+k_3x+k_4y)$. Section 4 contains a
few examples for data defined on oval and rounded corner domains.

\section{Equivalence transformations related to (\ref{eq2})}

Let us consider a one-parameter Lie group of transformations acting
on an open set ${\cal {D}}\subset \Omega \times \cal {W} \times
{\cal {E}}$, where ${\cal {W}}$ is the space of the data functions,
and ${\cal {E}}$ is the space of the parameter functions, given by
\begin{equation}\label{ET1}
\left\{\!\!
\begin{array}{l}
\tilde {x}=x+\varepsilon \Gamma (x,y,w,E)+{\cal O}(\varepsilon
^2),\\[6pt]
\tilde {y}=y+\varepsilon \Lambda (x,y,w,E) + {\cal O} (\varepsilon
^2),\\[6pt]
\tilde {w}=w+\varepsilon \Phi (x,y,w,E) +{\cal O} (\varepsilon
^2),\\[6pt]
\tilde {E}=E+\varepsilon \Psi (x,y,w,E) +{\cal O} (\varepsilon ^2),
\end{array}\right.
\end{equation}
where $\varepsilon $ is the group parameter. Let
\begin{equation}\label{vectorET1}
{\cal {V}}=\Gamma (x,y,w,E)\partial _x+\Lambda (x,y,w,E)\partial _y+
\Phi (x,y,w,E)
\partial _w +\Psi (x,y,w,E) \partial _E
\end{equation}
be its associated general infinitesimal generator. Assume that
$E=E(x,y)$ and $w=w(x,y)$ are both dependent variables in
(\ref{eq2}). The transformation (\ref{ET1}) is called an
\textit{equivalence transformation} related to the PDE (\ref{eq2})
if this leaves the equation invariant, i.e., $\tilde w_{\tilde x}
\tilde E_{\tilde x}+\tilde w_{\tilde y}\tilde E_{\tilde y}+\tilde
E(\tilde w_{\tilde x\tilde x}+\tilde w_{\tilde y\tilde y})=-1$. Note
that the set of all analytical solutions of (\ref{eq2}) will also be
invariant. For this case, the criterion for infinitesimal invariance
is given by
\[
\hbox{pr} ^{(2)}{\cal {V}}(F)\vert _{F=0}=0,
\]
with
\[
F\left(x,y,w^{(2)},E^{(1)}\right)=w_xE_x+w_yE_y+E(w_{xx}+w_{yy})+1,
\]
where $\hbox{pr} ^{(2)}{\cal {V}}$ denotes the second order
prolongation of the vector field ${\cal {V}}$ \cite{olverbook}.
Notice that this prolongation is determined by taking into account
that $E$ and $w$ are both dependent variables, exactly as one would
proceed in finding the classical Lie symmetries for a PDE without
arbitrary functions. The order of the prolongation of the vector
field ${\cal {V}}$ is given by the highest leading derivative of the
dependent variables. Applying the classical Lie method, we obtain
the following infinitesimals
\begin{equation}\label{ET2}
\left\{
\begin{array}{l}
\Gamma (x,y,w,E) = k_1-k_3y+k_4x \\
\Lambda (x,y,w,E)  = k_2+k_3x+k_4 y\\
\Phi (x,y,w,E)  =  \mu (w) \\
\Psi (x,y,w,E) =  E(2k_4-\mu'(w)),
\end{array}\right.
\end{equation}
\par
\noindent  where $k_i$, $i=1,...,4$ are real constants and $\mu =
\mu (w)$ is an arbitrary function. Hence, the infinitesimal
generator (\ref{vectorET1}) becomes
\begin{subequations}\label{vectorET2}
\begin{equation}
{\cal {V}}=\sum _{i=1} ^4 k_i {\cal {V}}_i+{\cal {V}}_{\mu },
\end{equation}
where
\begin{equation}
\begin{array}{c}
{\cal {V}}_1=\partial _x, \quad {\cal {V}}_2=\partial _y,\quad {\cal
{V}}_3=-y
\partial _x +x\partial _y,\quad {\cal {V}}_4=x\partial _x+y\partial
_y+2E\partial _E,\\[6pt]
{\cal {V}}_{\mu }=\mu (w) \partial _w -E \mu'(w)\partial _E.
\end{array}
\end{equation}
\end{subequations}
We obtain the following result

\begin{proposition}\label{prop}
There is an infinite dimensional Lie algebra of the equivalence
transformations related to (\ref{eq2}) spanned by the infinitesimal
generators\/ {\rm(\ref{vectorET2})}.
\end{proposition}

Therefore, the PDE {\rm(\ref{eq2})} is invariant under translations
in the $x$-space, $y$-space, rotations in the $(x,y)$-space, and,
respectively, scaling transformations in the\/ $(x,y,E)$-space. For
instance, we find that (\ref{eq2}) is invariant under translations
in $w$-space if we choose $\mu = const.$, and that the equation
remains unchanged under scaling transformations in $(w,E)$-space if
we choose $\mu(w) = w$.

Furthermore, an equivalence transformation (\ref{ET1}) can be used
to reduce the dimension of (\ref{eq2}) by augmenting it with
\[
\left\{\!\!
\begin{array}{l}
\Gamma (x,y,w,E)w_x+\Lambda (x,y,w,E)w_y- \Phi (x,y,w,E)=0,\\[6pt]
\Gamma (x,y,w,E)E_x+\Lambda (x,y,w,E)E_y - \Psi (x,y,w,E)=0,
\end{array}
\right.
\]
which is a first order PDE system defining the characteristics of
the vector field (\ref{vectorET1}). These relations are also called
\textit{invariance surface conditions}.

\section{Nonclassical equivalence transformations related to (\ref{eq2})}

Consider a one-parameter Lie group of transformations acting on an
open set ${\cal {D}}\subset \Omega \times \cal {W}\times {\cal
{E}}$, where $\cal {W}$ is the space of the data functions, and
$\cal {E}$ is the space of the parameter functions, given by
\begin{equation}\label{NET1}
\left\{\!\!
\begin{array}{l}
\tilde {x} =x+\varepsilon \xi(x,y,w,E)+{\cal O}(\varepsilon
^2),\\[6pt]
\tilde {y} =y+\varepsilon \eta (x,y,w,E) + {\cal O} (\varepsilon
^2),\\[6pt]
\tilde {w}=w+\varepsilon \phi (x,y,w,E) +{\cal O} (\varepsilon
^2),\\[6pt]
\tilde {E} =E+\varepsilon \psi (x,y,w,E) +{\cal O} (\varepsilon ^2),
\end{array}\right.
\end{equation}
where $\varepsilon $ is the group parameter. Let the following
vector field
\begin{equation}\label{vectorNET1}
{\cal {U}}=\xi (x,y,w,E)\partial _x+\eta (x,y,w,E)\partial _y+ \phi
(x,y,w,E)
\partial _w +\psi (x,y,w,E) \partial _E
\end{equation}
be the general infinitesimal generator related to (\ref{NET1}). The
transformation (\ref{NET1}) is called a \textit{nonclassical
equivalence transformation} (or a \textit{conditional symmetry}) of
the PDE (\ref{eq2}) if this leaves the subset
\[
S_{F,\phi _1,\phi _2}=\{F(x,y,w^{(2)},E^{(2)})=0,\\ \phi
_1(x,y,w^{(1)},E^{(1)})=0,\phi _2(x,y,w^{(1)},E^{(1)})=0\}
\]
of the set of all analytical solutions invariant, where
\begin{equation}
\label{ISCEW1}
\left\{\!\!
\begin{array}{l}
\phi _1:=\xi (x,y,w,E)w_x+\eta (x,y,w,E)w_y- \phi (x,y,w,E) =0,\\[6pt]
\phi _2:=\xi(x,y,w,E)E_x+\eta (x,y,w,E)E_y - \psi (x,y,w,E) =0
\end{array}
\right.
\end{equation}
represents the characteristics of the vector field ${\cal {U}}$ (or
the \textit{invariant surface conditions}). Here the criterion for
infinitesimal invariance is the following
\[
\left\{
\begin{array}{l}
\hbox{pr}^{(2)}{\cal {U}}(F)\vert_{F=0,\phi _1=0,\phi _2=0} =0,\\[6pt]
\hbox{pr}^{(1)}{\cal {U}}(\phi _1)\vert_{F=0,\phi _1=0,\phi _2=0} =0,\\[6pt]
\hbox{pr}^{(1)}{\cal {U}}(\phi _2)\vert_{F=0,\phi _1=0,\phi _2=0}
=0.
\end{array}
\right.
\]
If $\eta \neq 0 $, one can assume without loss of generality that
$\eta =1 $ (the case $\eta =0 $ is not discussed in this paper),
and, hence, (\ref{ISCEW1}) turns into
\begin{equation}
\label{ISCEW2}
\left\{
\begin{array}{l}
w_y =\phi (x,y,w,E)-\xi (x,y,w,E)w_x,\\
E_y=\psi (x,y,w,E)-\xi (x,y,w,E)E_x.
\end{array}
\right.
\end{equation}
At the first step, we augment the original PDE with (\ref{ISCEW2})
and eliminate all the partial derivatives of $w$ and $E$ with
respect to $y$ occurring in (\ref{eq2}). Hence, by using
(\ref{ISCEW2}) and its differential consequences, we obtain
\begin{subequations}\label{reduced1}
\begin{equation}
{\cal{A}}_1w_{xx}+{\cal{A}}_2w_x^2+{\cal{A}}_3w_x E_x+{\cal{A}}_4w_x+{\cal{A}}_5E_x+{\cal{A}}_6=0,
\end{equation}
where the coefficients ${\cal{A}}_i={\cal{A}}_i(x,y,w,E)$,
$i=1...6$, are the following
\begin{equation}
\begin{array}{l}
{\cal{A}}_1 =E(\xi ^2+1),\\
{\cal{A}}_2 =2E\xi \xi _w ,\\
{\cal{A}}_3 =\xi ^2+2E\xi \xi _E +1,\\
{\cal{A}}_4 =-\xi \psi +E\left ( \xi \xi _x - \xi _y -\phi \xi _w - \psi \xi _E-2\xi \phi _w\right ),\\
{\cal{A}}_5 = -\xi \left ( \phi +2E\phi _E \right ),\\
{\cal{A}}_6 = \phi \psi +1 -E\left (\xi \phi _x -\phi _y -\phi \phi
_w - \psi \phi _E\right ).
\end{array}
\end{equation}
\end{subequations}
The equation (\ref{reduced1}) has been obtained by using the
\texttt{GENDEFNC} command
\begin{verbatim}
   gendefnc(PDE, [w,E], [x,y], y, 3).
\end{verbatim}
Since ${\cal{A}}_1\neq 0$, (\ref{reduced1}) may be regarded as an
ODE in the unknown functions $w$ and $E$ (with $y$ as a parameter).
At the second step, by using the \texttt{GENDEFNC} command
\begin{verbatim}
   gendefnc(PDE, [w,E], [x,y], y)
\end{verbatim}
we obtain the \textit{determining equations} of the nonclassical
symmetries. This is an over-determined nonlinear PDE system for the
infinitesimals $\xi =\xi(x,y,w,E)$ and $\phi =\phi (x,y,w,E)$. Among
these equations, we get
\[
\xi _w=0, \quad \xi _E=0, \quad \phi _E=0,
\]
which implies
\begin{equation}
\label{xi1}
\xi (x,y,w,E) = A(x,y),\\[6pt]
\end{equation}
and
\begin{equation}
\label{phi1} \phi (x,y,w,E)  = G(x,y,w),
\end{equation}
where $A$ and $G$ are arbitrary functions. The substitution of the
above functions into the remaining equations yields
\begin{equation}
\label{psi1} \psi (x,y,w,E) = EF(x,y,w),
\end{equation}
with $F$ an arbitrary function of its arguments. By using the above
relations, the determining system is reduced to
\begin{equation}
\label{eqnF} F=-G_w+{2(A_x-AA_y)\over {A^2+1}},
\end{equation}
\begin{equation}
\label{eqnG1}
G_x-AG_y-{2AA_x-A_y(A^2-1)\over {A^2+1}}G=0,
\end{equation}
\begin{equation}
\label{eqnG2}
\begin{array}{l}
G_{xx}+G_{yy}+F_yG+FG_y+2GG_{yw}\\[6pt]
+{2(A_x-AA_y)\over {A^2+1}}(GG_w+G_y+FG)-{2(AA_x+A_y)\over
{A^2+1}}G_x=0,
\end{array}
\end{equation}
and
\begin{equation}\label{eqnA}
\begin{array}{l}
(A^4-1)A_{xx}+4A(A^2+1)A_{xy}-(A^4-1)A_{yy}-2A(A^2-3)A_x^2
\\[6pt]
-4(3A^2-1)A_xA_y+2A(A^2-3)A_y^2=0,
\end{array}
\end{equation}
where $A=A(x,y)$, $G=G(x,y,w)$, and $F=F(x,y,w)$ are the unknown
functions.

\medskip
\par
\noindent \textit{Remark}. The nonclassical symmetries do not leave
the form of the equation (\ref{eq2}) invariant. Moreover, the
nonclassical operators (\ref{vectorNET1}) do not form a vector
space, still less a Lie algebra, as the symmetry operators do. Since
every classical symmetry is a nonclassical symmetry but not
conversely, there exists a set of common solutions of the
determining system of the nonclassical symmetries (nonclassical
equivalence symmetries, in this context) and the determining system
of the classical symmetries (equivalence transformations,
respectively). This common solution is given by
\begin{equation}\label{comparison}
\left\{
\begin{array}{l}
\xi (x,y,w,E) = {\Gamma (x,y,w,E)\over {\Lambda
(x,y,w,E)}}={k_1-k_3y+k_4x\over {k_2+k_3x+k_4 y}}, \\[6pt]
\phi (x,y,w,E) = {\Phi (x,y,w,E)\over {\Lambda (x,y,w,E)}}= {\mu
(w)\over
{k_2+k_3x+k_4 y}},\\[6pt]
\psi (x,y,w,E )= {\Psi (x,y,w,E)\over {\Lambda
(x,y,w,E)}}={E(2k_4-\mu'(w))\over {k_2+k_3x+k_4
y}},\end{array}\right.
\end{equation}
where at least one of the constants $k_2,\;k_3$ or $k_4$ is nonzero.

To solve the determining equations (\ref{eqnF})--(\ref{eqnA}), we
proceed as follows: the equations (\ref{eqnG1}) and (\ref{eqnG2})
are analyzed in \S 3.1, the equation (\ref{eqnA}) is studied in \S
3.2, and the solutions of the determining equations
(\ref{eqnF})--(\ref{eqnA}) are given in \S 3.3.

\subsection{The equations (\ref{eqnG1}) and (\ref{eqnG2}) for $\phi =G(x,y,w)$}
First observe that (\ref{eqnG1}) can be written in the following
conservation form
\[
\left({G\over{A^2+1}}\right )_x - \left({AG\over {A^2+1}}\right)_y =
0.
\]
If $A\not \equiv 0$, then two cases might occur:
\medskip
\par
\noindent
\textit{Case G1.} $G\equiv 0$.  In this case, (\ref{eqnG1}) and (\ref{eqnG2}) are both satisfied.

\medskip
\par
\noindent \textit{Case G2.} $G \not \equiv 0$. There exists a
potential function $K=K(x,y,w)$ such that
\[
\left\{
\begin{array}{l}
K_x ={AG\over {A^2+1}}\\[10pt]
K_y={G\over {A^2+1}}. \end{array}\right.
\]
The above system yields $K_x=AK_y$ whose general solution is
\[
K(x,y,w)=P(u,w),
\]
where $P$ is an arbitrary function and $u=u(x,y)$ is a solution of
the equation
\begin{equation}
\label{eqU} u_x=A(x,y)u_y.
\end{equation}
From these relations we obtain
\begin{equation}
\label{solGnew}
G=u_y(A^2+1)S,
\end{equation}
where $S(u,w)$ denotes the partial derivative $P_u(u,w)$. Thus, the
solution of (\ref{eqnG1}) is given by (\ref{solGnew}), where $u$
satisfies (\ref{eqU}).

The substitution of (\ref{solGnew}) into (\ref{eqnG2}) implies
\begin{equation}\label{eqnS}
q_1(S_{uw}S-S_uS_w+S_{uu})+q_2S_u+q_3S=0,
\end{equation}
where the coefficients $q_i$ are expressed in terms of $A$ and $u_y$
as follows
\[
\begin{array}{l}
q_1 =u^3_y(A^2+1)^3,\\[6pt]
q_2 =u_y(A^2+1)\left[3u_{yy}(A^2+1)^2+u_y(5A_x+3A^3A_y+3A^2A_x+
AA_y)\right],\\[6pt]
q_3 =u_{yyy}(A^2+1)^3+u_{yy}(A^2+1)[(3A^2+5)A_x+(3A^2+1)AA_y]\\[6pt]
\;\;\;\;\;\;+u_y[2A(A^2+1)A_{xx}+(A^2+1)(A^2+3)A_{xy}+A(A^2+1)^2A_{yy}\\[6pt]
\;\;\;\;\;\;-2(A^2-3)A_x^2+2A(A^2-3)A_xA_y+(A^4-1)A_y^2].
\end{array}
\]
In particular, the method of separation of variables applied to
(\ref{eqnS}) implies solutions of the form $S(u,w) = p(u)\mu (w)$,
where $p=p(u)$ satisfies the equation
\begin{equation}\label{eqnp}
q_1p_{uu}+q_2p_u+q_3p=0
\end{equation}
and $\mu =\mu (w)$ is an arbitrary function of its argument.

To summarize, the equation (\ref{eqnG2}) has been reduced to
(\ref{eqnS}), where $A=A(x,y)$ satisfies (\ref{eqnA}) and $u=u(x,y)$
is a solution of (\ref{eqU}).

\subsection{The equation (\ref{eqnA}) for $\xi =A(x,y)$}
In the following, we show that (\ref{eqnA}) can be reduced to a
Monge-Amp\`{e}re equation.

\medskip
\par
\noindent \textit{Case A1}. $A\equiv 0$
\medskip
\par
\noindent The equation (\ref{eqnA}) admits the trivial solution
$A\equiv 0$.
\medskip
\par
\noindent \textit{Case A2}. $A=k$, with $k\neq 0$
\medskip
\par
\noindent The constant function $A=k$, with $k\neq 0$ is also a
solution of (\ref{eqnA}).
\medskip
\par
\noindent \textit{Case A3}. $A$ is a nonconstant function
\medskip
\par
\noindent If $A$ is a nonconstant function, then an alternative
formulation of (\ref{eqnA}) is
\[
\begin{array}{cc}
B(B^2+1)B_{xx}+(B^2-1)(B^2+1)B_{xy}-B(B^2+1)B_{yy}\\[10pt]
\quad -(3B^2-1)B_x^2-2B(B^2-3)B_xB_y+(3B^2-1)B_y^2=0
\end{array}
\]
where $B=(A+1)/(A-1)$. The conservation form of the above PDE is
\begin{equation}
\label{eqforBcons} \left( {B(B_x+BB_y)\over {(B^2+1)^2}}\right )_x
-\left ( {B_x+BB_y\over {(B^2+1)^2}}\right )_y=0.
\end{equation}
We distinguish the following two cases:
\medskip
\par
\noindent \textit{Case A3.1.} $B=B(x,y)$ is a solution of
(\ref{eqnMonge}).
\medskip
\par
\noindent If $B=B(x,y)$ satisfies the Monge equation
\begin{equation}
\label{eqnMonge} B_x+BB_y=0,
\end{equation}
then (\ref{eqforBcons}) holds. Since the general solution of
(\ref{eqnMonge}) is given implicitly by
\[
y-xB=\nu(B),
\]
where $\nu $ is an arbitrary function, the corresponding solution of
(\ref{eqnA}) is
\begin{equation}\label{Anew}
y-x{A+1\over {A-1}}=\nu \left({A+1\over {A-1}}\right).
\end{equation}

\medskip
\par
\noindent \textit{Case A3.2.} $B=B(x,y)$ does not satisfy
(\ref{eqnMonge})
\medskip
\par
\noindent It follows from (\ref{eqforBcons}) that there exists a
potential function $T=T(x,y)$ such that
\begin{equation} \label{equationsT}
\left\{
\begin{array}{l}
T_x={B_x+BB_y\over {(B^2+1)^2}},\\[10pt]
T_y={B(B_x+BB_y)\over {(B^2+1)^2}}.
\end{array}\right.
\end{equation}
Since $T_x$ cannot be identically zero, the above equations yield
$B=T_y/T_x$ and, by substituting it into the first equation of
(\ref{equationsT}) we have
\[
T_xT_yT_{xx}-\left (T_x^2-T_y^2\right )T_{xy}-T_xT_yT_{yy}+\left (T_x^2+T_y^2 \right ) ^2=0.
\]
By using the following Legendre transformation
\cite[p.~\mbox{353}]{zwillinger}
\[
{\cal {H}}(a,b)+T(x,y)=xa+yb,
\]
where $T_x=a$, $x={\cal {H}}_a$, $T_y=b$, and $y={\cal {H}}_b$, the
above PDE turns into the following Monge-Amp\`{e}re equation
\[
{\cal {H}}_{aa}{\cal {H}}_{bb}-{\cal {H}}_{ab}^2-{ab\over {\left
(a^2+b^2\right )^2}}{\cal {H}}_{aa}+{a^2-b^2\over {\left
(a^2+b^2\right )^2}}{\cal {H}}_{ab}+ {ab \over {\left (a^2+b^2\right
)^2}}{\cal {H}}_{bb}=0.
\]
Furthermore, this can be reduced to the Monge-Amp\`{e}re equation
\begin{equation}\label{MAinV}
V_{aa}V_{bb}-V_{ab}^2=-{1\over {\left (a^2+b^2\right )^2}},
\end{equation}
where
\begin{equation}\label{changeH}
V(a,b)={\cal {H}}(a,b)-{1\over {2}}\arctan \left ({a\over {b}}\right
).
\end{equation}
Exact solutions for particular Monge-Amp\`{e}re equations have been
extensively analyzed in \cite{handbookPZ2}. The Monge-Amp\`{e}re
equation (\ref{MAinV}) may be included in Case 17
\cite[p.~\mbox{458}]{handbookPZ2} or Case 20
\cite[p.~\mbox{460}]{handbookPZ2}.

\subsection{Solutions of the determining equations (\ref{eqnF})--(\ref{eqnA})}
We distinguish the following four cases:

\medskip
\par
\noindent \textit{Case 1}. $A\equiv 0$ and $G \equiv 0$
\medskip
\par
\noindent

\begin{proposition}\label{prop1}
If $A\equiv 0$ and $G \equiv 0$, the infinitesimal generator
(\ref{vectorNET1}) becomes ${\cal {U}}=\partial _y$, which implies
the invariance of the equation (\ref{eq2}) with respect to
translations in the $y$-space.
\end{proposition}

This nonclassical symmetry reduction is an equivalence
transformation that can be obtained from (\ref{comparison}) for $\mu
=0$, $k_2=1$, and $k_i=0$, where $i=1,3,4$.

\medskip
\par
\noindent \textit{Case 2}. $A\equiv 0$ and $G \not \equiv 0$
\medskip
\par
\noindent The PDE (\ref{eqnG1}) takes the form $G_x=0$ and it
follows that $G=H(y,w)$. After substituting it into (\ref{eqnF}) and
(\ref{eqnG2}), we get $F=-H_w$ and, respectively,
\begin{equation}\label{eqnGyw}
H_{yy}+HH_{yw}-H_yH_w=0,
\end{equation}
where $H$ is an arbitrary function. Notice that $H\not \equiv 0$ in
the above equation. We distinguish the following cases:
\medskip
\par
\noindent \textit{Case 2.1}. $H=\mu (w)$, where $\mu $ is an
arbitrary function
\medskip
\par
\noindent  Indeed, $H=\mu (w)$ is a particular solution of
(\ref{eqnGyw}).

\begin{proposition}\label{prop20}
If $A\equiv 0$ and $G=\mu (w)$, where $\mu $ is an arbitrary
function, the infinitesimal generator (\ref{vectorNET1}) turns into
\begin{equation}
\label{another} {\cal {U}}=\partial _y+\mu (w)\partial _w-\mu
'(w)E\partial _E.
\end{equation}
\end{proposition}

The above nonclassical equivalence transformation is in fact an
equivalence transformation and corresponds to the case $k_2=1$, and
$k_i=0$, where $i=1,3,4$ in (\ref{comparison}).
\medskip
\par
\noindent \textit{Case 2.2}. $H_y\not \equiv 0$
\medskip
\par
\noindent The PDE (\ref{eqnGyw}) can be written in the following
conservation form
\[
\left({1\over {H_y}}\right )_y + \left ( {H\over {H_y}} \right)_w =0.
\]
After introducing the potential function $g=g(y,w)$, we get
\[
\left\{
\begin{array}{l}
g_y = {H\over {H_y}} \\[6pt]
g_w = -{1\over {H_y}}.
\end{array}\right.
\]
Clearly, $g_w\not \equiv 0$. After eliminating $H_y$ in the above
system, we obtain $H=-g_y/ g_w$. Next, substituting it into the
second equation, the following PDE results
\begin{equation}\label{eqnforg}
g_y(g_y-y)_w-g_w(g_y-y)_y=0.
\end{equation}
The following two cases occur:

\medskip
\par
\noindent \textit{Case 2.2.1}. $g_y=y$
\medskip
\par
\noindent It results $g(y,w)=y^2/2+h(w)$, where $h$ is an arbitrary
nonconstant function (otherwise, $g_w \equiv 0$). Since $H=-g_y/
g_w$, we get $H=y\mu (w)$, where $\mu =-1/h'$.
\begin{proposition}\label{prop3}
If $A\equiv 0$ and $G=y\mu (w)$, where $\mu $ is an arbitrary
function, then the infinitesimal generator (\ref{vectorNET1})
becomes
\begin{equation}\label{case221}
{\cal {U}}=\partial _y+ y \mu (w) \partial _w -y\mu '(w) E \partial
_E,
\end{equation}
\end{proposition}

The nonclassical symmetry generated by (\ref{case221}) is a new
symmetry reduction for (\ref{eq2}) which cannot be obtained from
(\ref{comparison}).

\medskip
\par
\noindent \textit{Case 2.2.2}. $g_y-y\not \equiv 0$
\medskip
\par
\noindent Since the equation (\ref{eqnforg}) is the Jacobian of the
functions $g$ and $g_y-y$, there exists a function $\alpha $ such
that
\begin{equation}\label{equation4g}
g_y=y+\alpha (g).
\end{equation}
In the above relation, $w$ is viewed as a parameter. The equation
(\ref{equation4g}) can be written as the following Abel ODE of
second kind
\[
\left (y+\alpha (g)\right
){dy\over {dg}}=1,
\]
for $y=y(g)$. The canonical substitutions $z=\alpha (g)$ and
$v=y+\alpha (g)$ reduce the above ODE to its canonical form
\begin{equation}\label{Abel10}
vv'-v=\beta (z),
\end{equation}
where $v=v(z)$ and $\beta =1/(\alpha ' \circ \alpha ^{-1})$. A
collection of the known cases of solvable Abel ODEs of the form
(\ref{Abel2}) is presented in \cite[pp.~\mbox{107--120}]{handbookPZ}
and new results can be found, for example, in \cite{Cheb1}. Each of
these ODEs corresponds to a nonclassical equivalence transformation
related to (\ref{eq2}).

\begin{proposition}\label{prop4}
If $A\equiv 0$ and $G=-g_y/g_w$, where $g$ is a solution of
(\ref{equation4g}), then the infinitesimal generator
(\ref{vectorNET1}) turns into
\begin{equation}\label{case222}
{\cal {U}}=\partial _y -g_y/g_w\partial _w + E(g_y/g_w)_w \partial
_E,
\end{equation}
\end{proposition}

Since $g$ satisfies (\ref{eqnforg}), the above vector field
generates new symmetry reductions related to (\ref{eq2}) that are
not equivalence transformations.
\medskip
\par
\noindent \textit{Case 3}. $A\not \equiv 0$ and $G \equiv 0$
\medskip
\par
\noindent By (\ref{eqnF}), we obtain $F=2(A_x-AA_y)/(A^2+1)$.

\begin{proposition}\label{prop5}
If $A\not \equiv 0$ is a solution of the equation (\ref{eqnA}) and
$G \equiv 0$, then the nonclassical infinitesimal generator
(\ref{vectorNET1}) becomes
\begin{equation}\label{case3}
{\cal {U}}=A(x,y)\partial _x +\partial _y + {2(A_x-AA_y)\over
{A^2+1}}E\partial _E.
\end{equation}
\end{proposition}

\medskip
\par
\noindent \textit{Case 3.1}. $A=k$, where $k\neq 0$
\medskip
\par
\noindent
\begin{proposition}\label{prop6}
If $A=k$, where $k\neq 0$ is a constant, the infinitesimal generator
(\ref{case3}) rewrites as ${\cal {U}}=k\partial _x +\partial _y$.
\end{proposition}

Replacing $\mu =0$, $k_2=1$, and $k_i=0$ with $i=1,3,4$ in
(\ref{comparison}), we obtain the symmetry reduction generated by
the above vector field, and, therefore, this is an equivalence
transformation.

\medskip
\par
\noindent \textit{Case 3.2}. $A$ is a nonconstant function
\medskip
\par
\noindent If $A$ is a nonconstant function, then we distinguish two
subcases:
\medskip
\par
\noindent
\textit{Case 3.2.1}. $A$ is given implicitly by (\ref{Anew}).
\medskip
\par
\noindent
\textit{Case 3.2.2}. The equation for $A$ is reduced to the Monge-Amp\`{e}re equation (\ref{MAinV}).
\medskip

In the above two cases, the nonclassical operator (\ref{case3})
generates new symmetry reductions for (\ref{eq2}) excepting the case
when $A(x,y)=(k_1-k_3y+k_4x)/(k_2+k_3x+k_4y)$.
\medskip
\par
\noindent \textit{Case 4}. $A\not \equiv 0$ and $G\not \equiv 0$

\medskip
\par
\noindent \textit{Case 4.1}. $A=k$, where $k\neq 0$
\medskip
\par
\noindent Without loss of generality, we consider the particular
solution $u(x,y)=kx+y$ of the equation (\ref{eqU}). Since $u_y=1$,
the relation (\ref{solGnew}) yields $G=(k^2+1)S$, where $S$ is a
nontrivial solution of the equation (\ref{eqnS}) which rewrites as
\begin{equation}\label{eqnSuw}
S_{uu}+SS_{uw}-S_uS_w=0.
\end{equation}
The following cases might occur:
\medskip
\par
\noindent \textit{Case 4.1.1}. $S=\mu (w)$, where $\mu $ is an
arbitrary function
\medskip
\par
\noindent In this case, (\ref{eqnSuw}) is satisfied and (\ref{eqnF})
implies $F=-\mu ' (w)$.

\begin{proposition}\label{prop411}
If $A=k$ with $k\neq 0$ and $G=(k^2+1)\mu (w)$, where $\mu $ is an
arbitrary function, then the nonclassical infinitesimal generator
(\ref{vectorNET1}) becomes
\begin{equation}\label{eqncase3}
{\cal {U}}=k\partial _x +\partial _y +(k^2+1)\mu (w)\partial _w
-(k^2+1)\mu '(w)E\partial _E.
\end{equation}
\end{proposition}

The above nonclassical operator generates an equivalence
transformation that can be obtained from (\ref{comparison}) for
$k_1=k/(k^2+1)$, $k_2=1/(k^2+1)$, $k_3=0$, and $k_4=0$.

\medskip
\par
\noindent \textit{Case 4.1.2}. $S_u\not \equiv 0$
\medskip
\par
\noindent Since the PDE (\ref{eqnSuw}) can be written in the
conservation form
\[
\left({1\over {S_u}}\right )_u + \left ( {S\over {S_u}} \right)_w
=0,
\]
there exists a potential function $Q=Q(u,w)$ such that
\[
\left\{
\begin{array}{l}
Q_u= {S\over {S_u}}, \\[6pt]
Q_w= -{1\over {S_u}}.
\end{array}\right.
\]
Eliminating $S$ in the above system, we obtain $S=-Q_u/Q_w$, where
\begin{equation}\label{eq4Q}
Q_u(Q_u-u)_w-Q_w(Q_u-u)_u=0.
\end{equation}
We distinguish the following two cases:
\medskip
\par
\noindent \textit{Case 4.1.2.a}. $Q_u=u$
\medskip
\par
\noindent In this case, $Q(u,w)=u^2/2+p(w)$, where $p$ is an
arbitrary nonconstant function (otherwise $Q_w\equiv 0$). With the
aid of $S=-Q_u/Q_w$, we get $S(u,w)=u\mu (w)$, where $\mu =-1/p'$
and $u(x,y)=kx+y$.

\begin{proposition}\label{412a}
For $A(x,y)=k$ ($k\neq 0$) and $G(x,y,w)=(k^2+1)(kx+y)\mu (w)$,
where $\mu $ is an arbitrary nonconstant function, the nonclassical
infinitesimal generator (\ref{vectorNET1}) rewrites as follows
\begin{equation}\label{prop412a}
{\cal {U}}=k\partial _x +\partial _y +(k^2+1)(kx+y)\mu (w)\partial
_w -(k^2+1)(kx+y)\mu '(w)E\partial _E.
\end{equation}
\end{proposition}

The nonclassical operator (\ref{prop412a}) generates new
transformations that are not equivalence transformations related to
(\ref{eq2}).

\medskip
\par
\noindent \textit{Case 4.1.2.b}. $Q_u-u\not \equiv 0$
\medskip
\par
\noindent Observe that (\ref{eq4Q}) is the Jacobian of $Q$ and
$Q_u-u$. Thus, there exists a function $\gamma $ such that
\begin{equation}
\label{equationg} Q_u=u+\gamma (Q).
\end{equation}
In the above ODE, $w$ is viewed as a parameter. Similar to Case
2.2.2, the equation (\ref{equationg}) can be reduced to
\[
\left (u+\gamma (Q)\right ){du\over {dQ}}=1,
\]
which is an Abel ODE of second kind for $u=u(Q)$. Moreover, after
using the substitutions $s=\gamma (Q)$ and $V=u+\gamma (Q)$, the
above ODE can be reduced to the canonical form
\begin{equation}\label{Abel2}
VV'-V=\theta (s),
\end{equation}
where $\theta =1/(\gamma ' \circ \gamma ^{-1})$. This is an Abel
equation of second kind for the unknown function $V=V(s)$. For each
solution of (\ref{Abel2}), we obtain a nonclassical symmetry for
(\ref{eq2}). The solvable Abel ODEs of the form (\ref{Abel2}) that
are known so far are listed in
\cite[pp.~\mbox{107--120}]{handbookPZ}. More recent results can be
found, for instance, in \cite{Cheb1}.
\par
\noindent
\begin{proposition}\label{prop412b}
If $A(x,y)=k$ with $k\neq 0$ and $G(x,y,w)=-(k^2+1)Q_u/Q_w$, where
$u(x,y)=kx+y$, and $Q=Q(u,w)$ satisfies (\ref{eq4Q}), then the
nonclassical infinitesimal generator (\ref{vectorNET1}) is given by
\begin{equation}\label{eqncase30}
{\cal {U}}=k\partial _x +\partial _y -(k^2+1)Q_u/Q_w\partial _w
+(k^2+1)E(Q_u/Q_w)_w\partial _E.
\end{equation}
\end{proposition}

Since $Q$ satisfies (\ref{eq4Q}), the above vector field generates
new symmetry reductions related to (\ref{eq2}) that are not
equivalence transformations.
\medskip
\par
\noindent \textit{Case 4.2}. Suppose $A$ is a nonconstant function
and $G\not \equiv 0$. \medskip
\par
\noindent
\begin{proposition}
If $A$ is a nonconstant function satisfying the equation
(\ref{eqnA}) and $G=u_y(A^2+1)S$, where $u$ satisfies (\ref{eqU})
and $S$ is a nonzero solution of (\ref{eqnS}), the nonclassical
infinitesimal generator (\ref{vectorNET1}) is written as
\begin{equation}\label{eqncase4}
{\cal {U}}=A\partial _x +\partial _y +u_y(A^2+1)S\partial _w
+E\left(-u_y(A^2+1)S_w+{2(A_x-AA_y)\over{A^2+1}}\right)\partial _E.
\end{equation}
\end{proposition}
Two cases might occur:
\medskip
\par
\noindent
\textit{Case 4.2.1}.  $A$ is given implicitly by (\ref{Anew}).

\medskip
\par
\noindent
\textit{Case 4.2.2}. The equation for $A$ is reduced to the Monge-Amp\`{e}re equation (\ref{MAinV}).
\medskip

In the above two cases, the vector (\ref{eqncase4}) generates new
symmetry reductions that cannot be obtained from (\ref{comparison})
excepting the case when $A(x,y)=(k_1-k_3y+k_4x)/(k_2+k_3x+k_4y)$.

We thus have found new symmetry reductions related to (\ref{eq2})
and have shown that these are given by (\ref{case221}),
(\ref{case222}), (\ref{prop412a}), (\ref{eqncase30}), and,
respectively, (\ref{case3}) and (\ref{eqncase4}) excepting the case
when $A(x,y)=(k_1-k_3y+k_4x)/(k_2+k_3x+k_4y)$.

\section{Examples}
The nonclassical symmetries related to (\ref{eq2}) yield classes of
data and suitable domains for which the dimension of the problem can
be reduced. Since the homogeneous Dirichlet boundary condition
(\ref{cond1}) is imposed, the data $w$ must satisfy (\ref{ISCEW2})
and the boundary $\partial \Omega $ must be "compatible" with the
symmetry reduction as well. To illustrate this approach, we discuss
a few examples of data modeled by the functions
\begin{equation}\label{data1}
w(x,y)=-y^2+W(x)
\end{equation}
which are invariant with respect to the nonclassical infinitesimal
generator
\begin{equation}
\label{newUk} {{\cal {U}}}=\partial _y -2y\partial _w,
\end{equation}
which is obtained from (3.27) by setting $\mu(w) = -2$. Substituting
(\ref{newUk}) in (\ref{ISCEW2}), we get $w_y=-2y$ and $E_y=0$. This
implies $E(x,y)= E_0(x)$ (i.e., the heat conductivity is constant
along the lines $x=const.$) and $w$ is given by (\ref{data1}).
Therefore, we can consider circles, ellipses, generalized Lam\'{e}
curves $x^{2p}+y^2=1$ ($p>2$), Granville's egg curve
$y^2x^2=(x-b)(1-x)$, where $b\neq 0,1$, or elliptical curves -- in
particular, Newton's egg curve $y^2=(x^2-1)(x-a)$, where $a\neq \pm
1$. In this case, as it is shown below, the parameter cannot be
determined at the points $(x_0,y)$ for which $W'(x_0)=0$. Despite
this fact, additional information about the parameter on the
specified domain can be obtained.

\begin{figure}[ht!]
\vspace{0.5cm} \epsfig{file=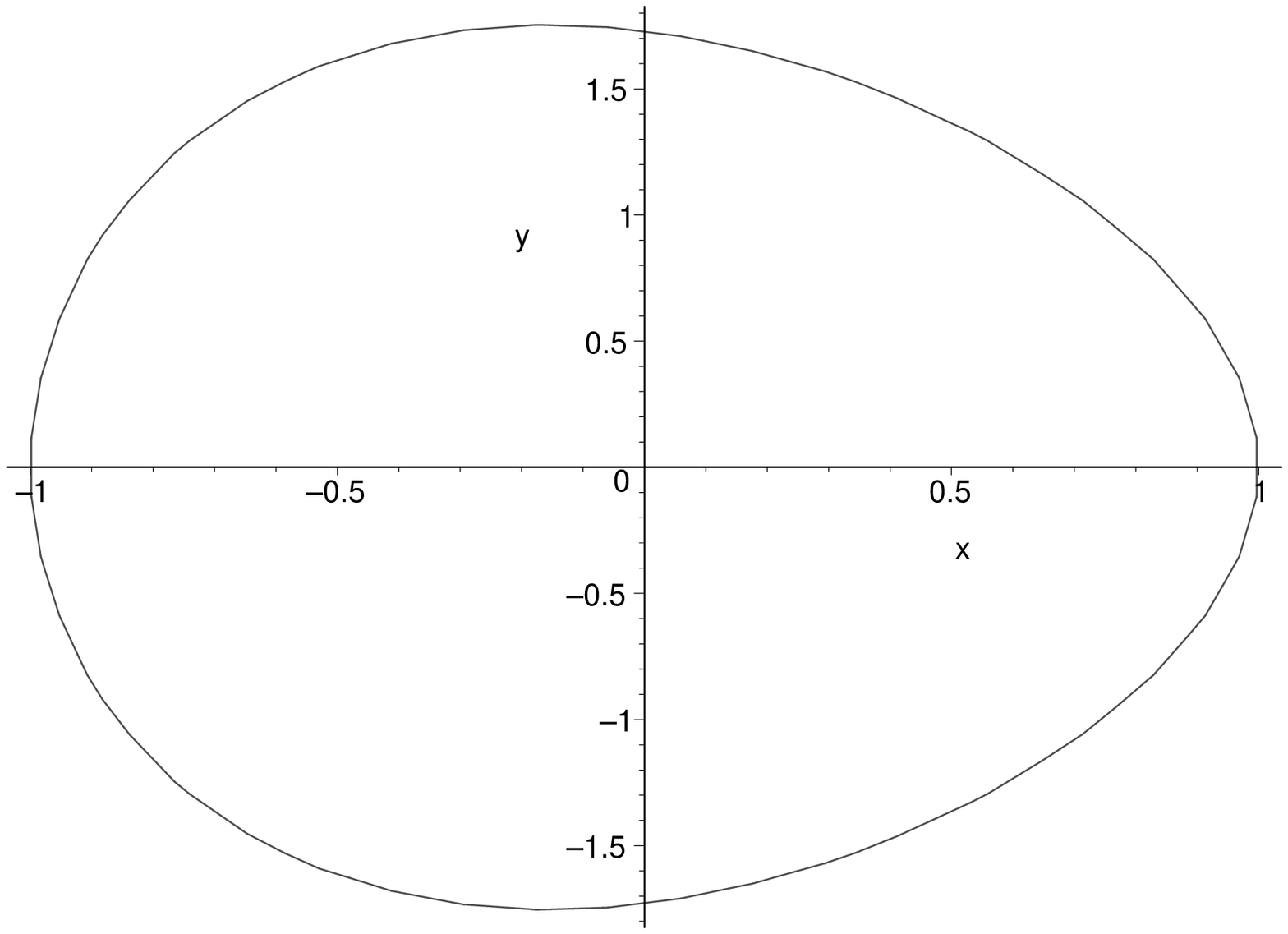,width=4cm, height=6cm,
angle=270} \hspace{0.75cm}
\epsfig{file=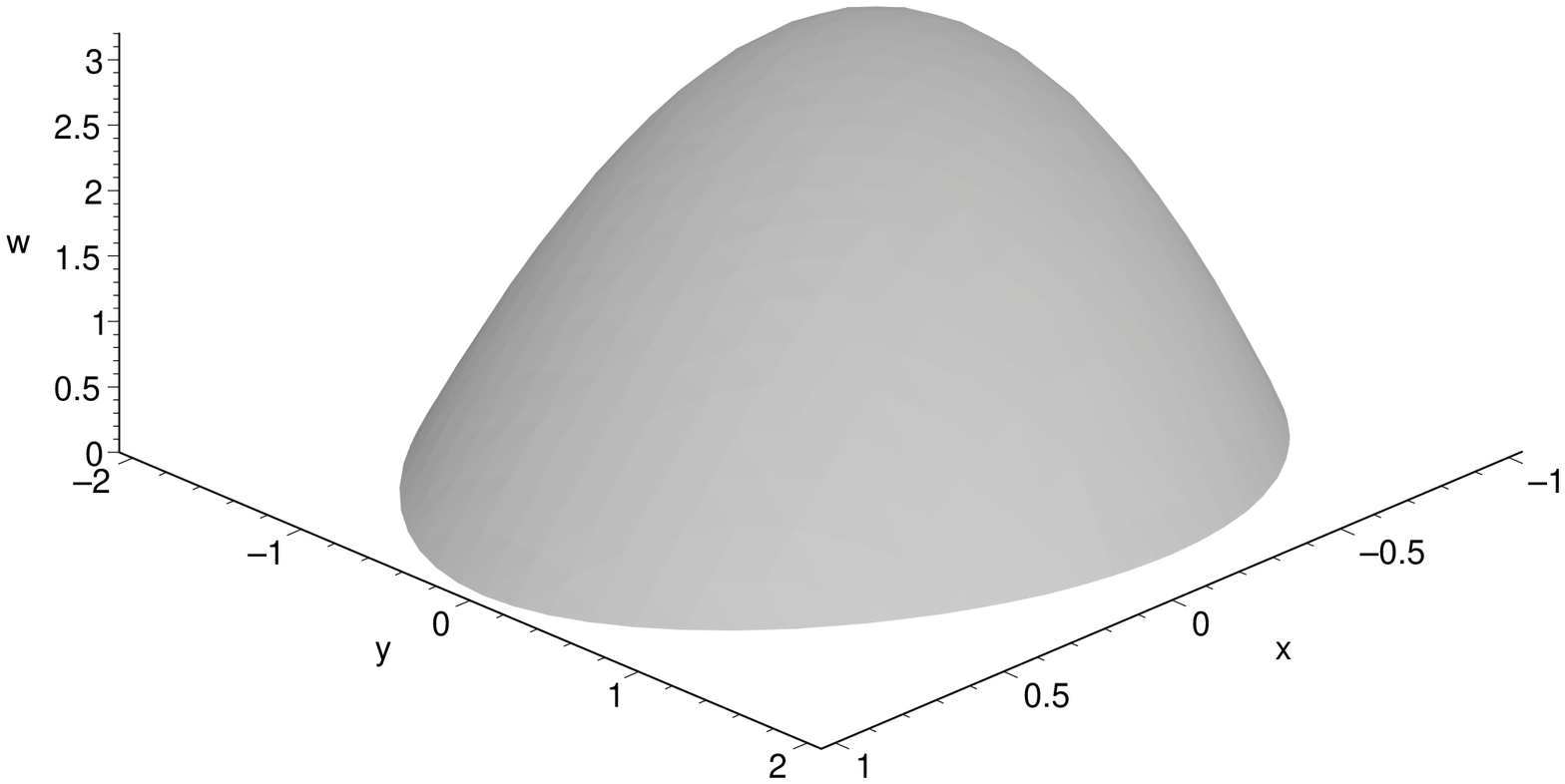,width=4cm,height=6cm,angle=-90} \vspace{1cm}
\caption{The boundary $\partial \Omega =\{(x,y):y^2-(x^2-1)(x-3)
=0,\;x\in [-1,1]\}$ and the data $w(x,y)=-y^2+(x^2-1)(x-3)$
discussed in Example 1.} \label{figure1}
\end{figure}

\begin{figure}[t!]
\vspace{0.5cm} \epsfig{file=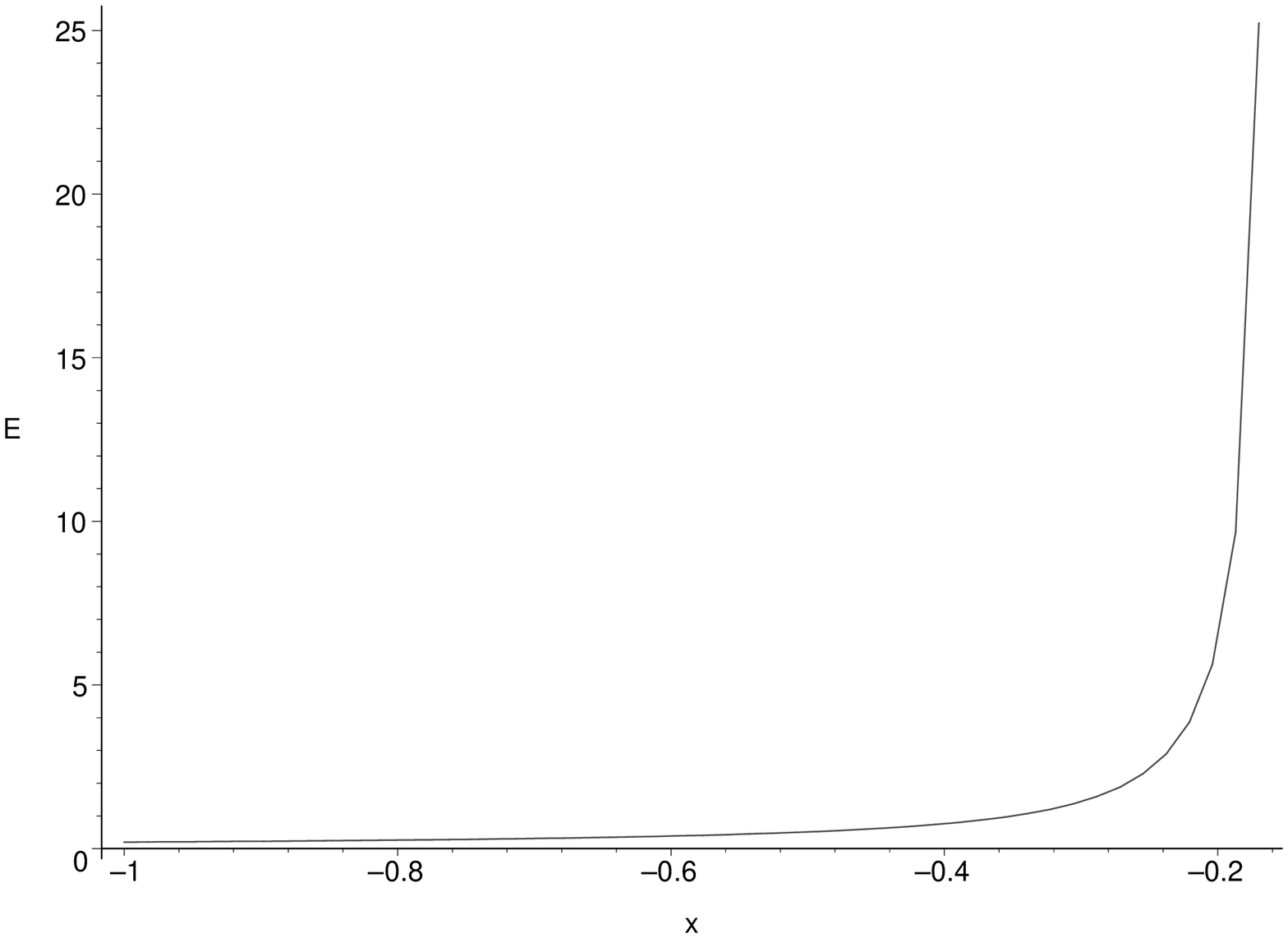,width=5cm, height=5cm,
angle=270} \hspace{1cm} \epsfig{file=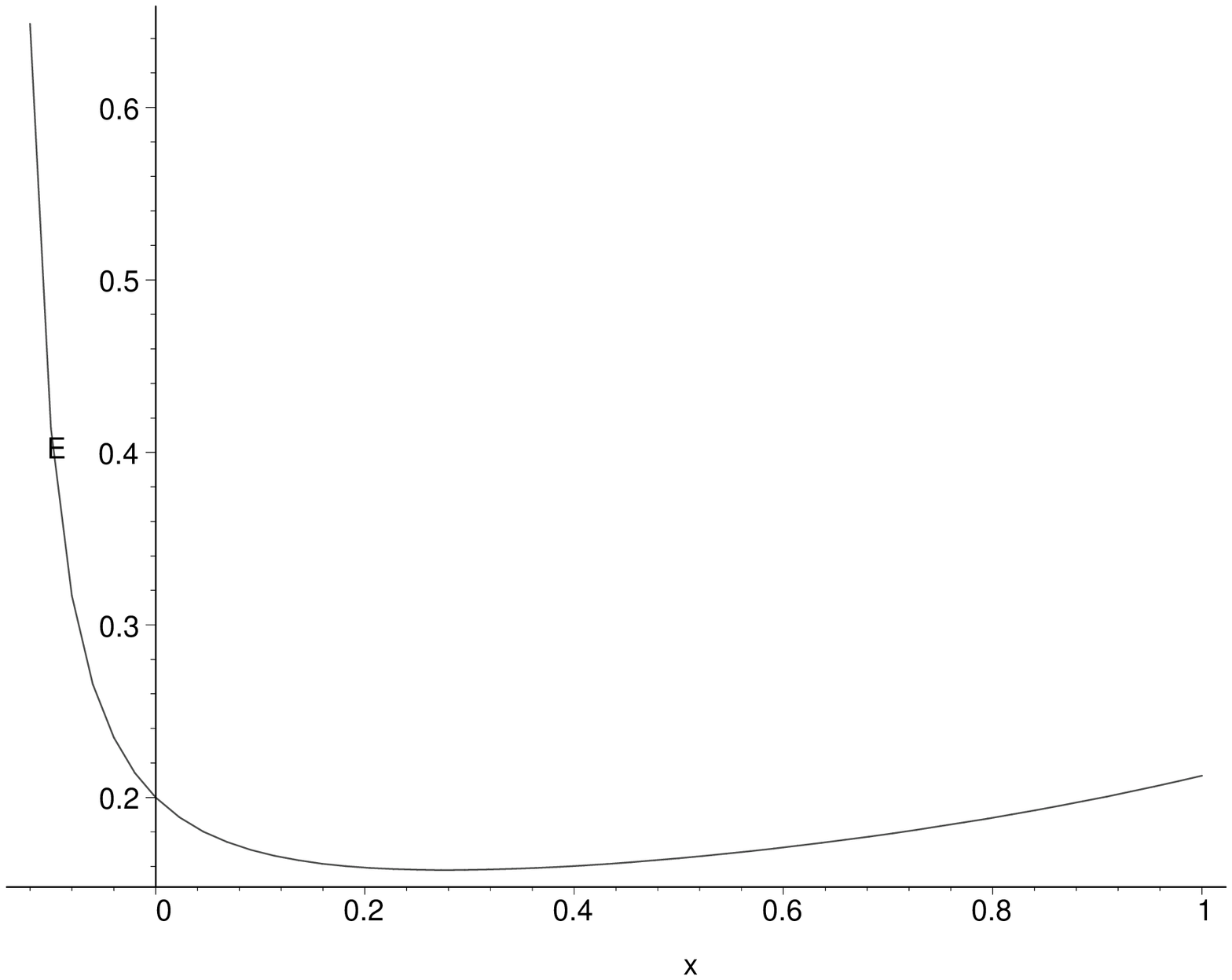,
height=4cm,width=5cm,angle=-90} \vspace{1cm} \caption{The parameter
$E(x,y)=E_0(x)$ corresponding to data discussed in Example 1.}
\label{figure2}
\end{figure}

\textit{Example 1}. (Newton's egg curve) Suppose $\partial \Omega
=\{(x,y):y^2=(x^2-1)(x-3),\;x\in [-1,1]\}$ (see Figure
\ref{figure1}). For the data
\[
w(x,y)=-y^2+(x^2-1)(x-3),
\]
the equation (\ref{eq2}) is reduced to the ODE
$(3x^2-6x-1)E_0'(x)+(6x-8)E_0(x)=-1$, whose general solution is
\[
E_0(x)={C-E_1(x)\over {3x^2-6x-1}}\exp \left [-{\sqrt{3}\over
{3}} \hbox{arctanh}\left({\sqrt{3}\over {2}}(x-1) \right)\right],
\]
for $x\neq 1-2/\sqrt{3}$, where $E_1'(x)=\exp \left [\sqrt{3}/3
\;\hbox{arctanh}\left((x-1)\sqrt{3}/2 \right)\right]. $ Figure
\ref{figure2} shows the parameters satisfying the conditions
$E(-1)=0.2$, and, respectively, $E(1)=0.2$. The parameter $E$ cannot
be determined for $x_0=1-2/\sqrt{3}$.

\vspace{0.5cm}

\begin{figure}[ht!]
\vspace{0.5cm} \epsfig{file=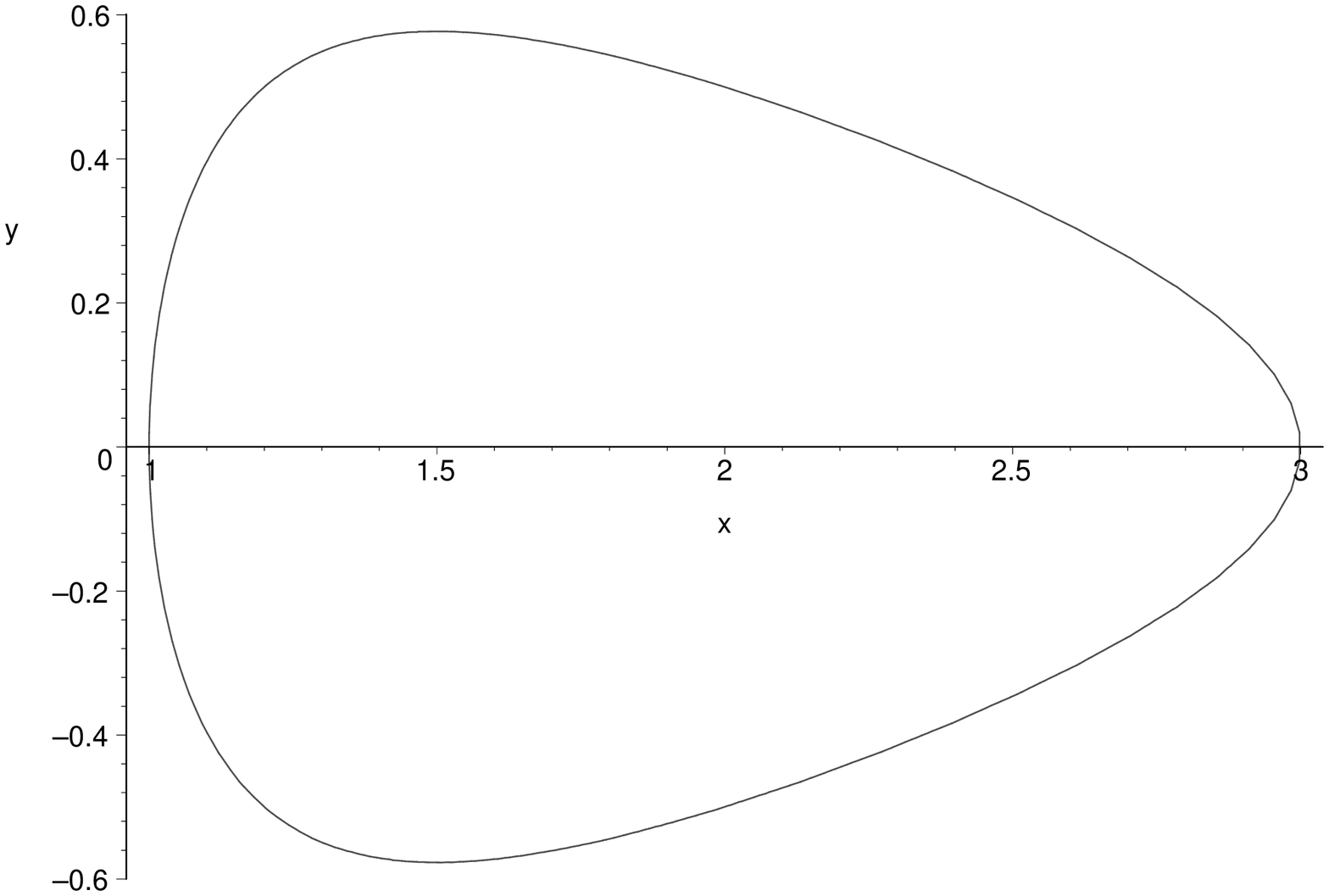,width=4cm, height=6cm,
angle=270} \hspace{0.75cm}
\epsfig{file=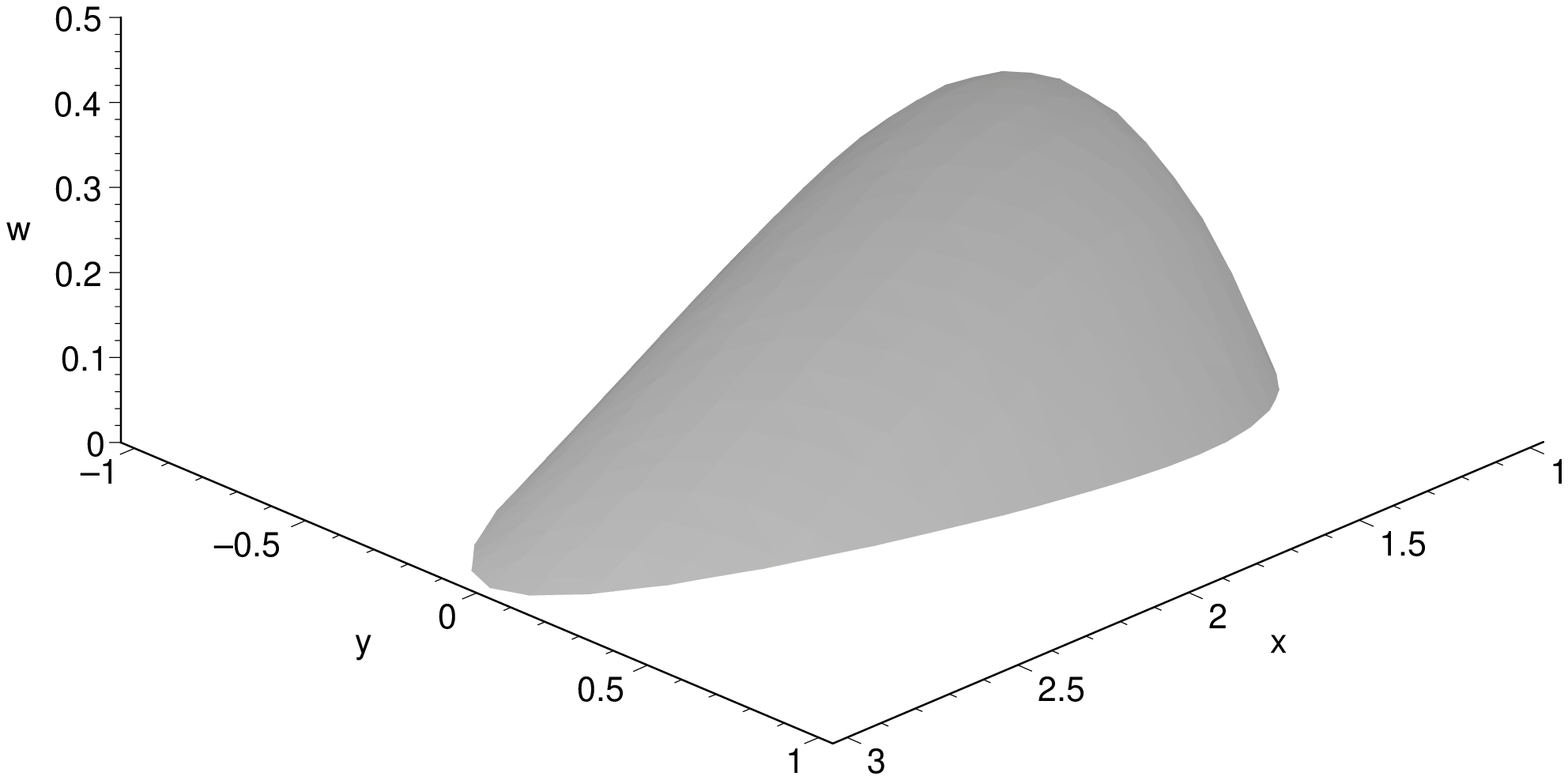,width=4cm,height=6cm,angle=-90} \vspace{1cm}
\caption{The boundary $\partial \Omega =\{(x,y):x^2y^2=(x-3)(1-x)
,\;x\in [1,3]\}$ and the data $w(x,y)=y^2-{(x-3)(1-x)\over {x^2}}$
considered in Example 2.} \label{figure3}
\end{figure}

\begin{figure}[t!]
\vspace{0.5cm} \epsfig{file=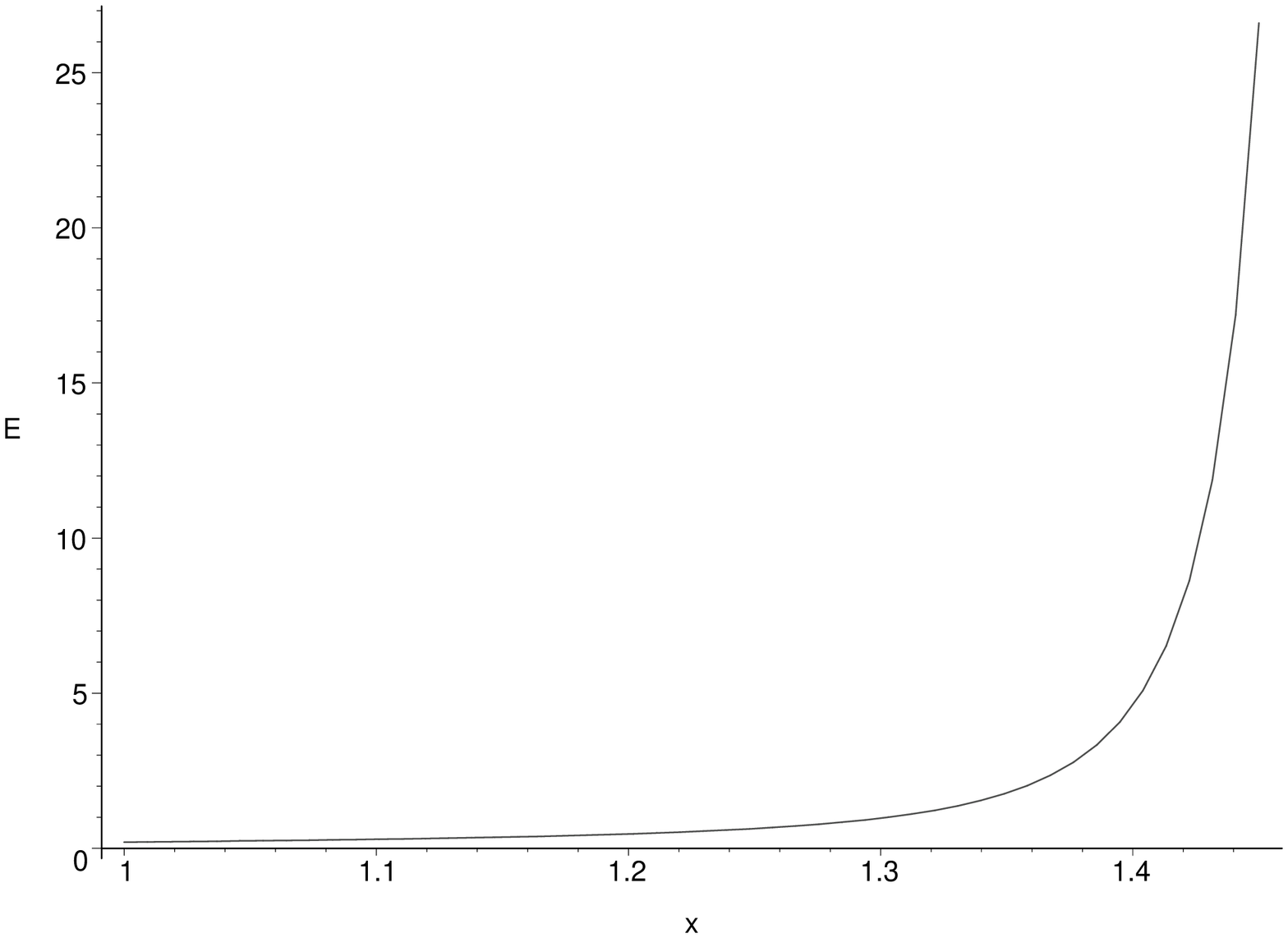,width=5cm, height=5cm,
angle=270} \hspace{1cm} \epsfig{file=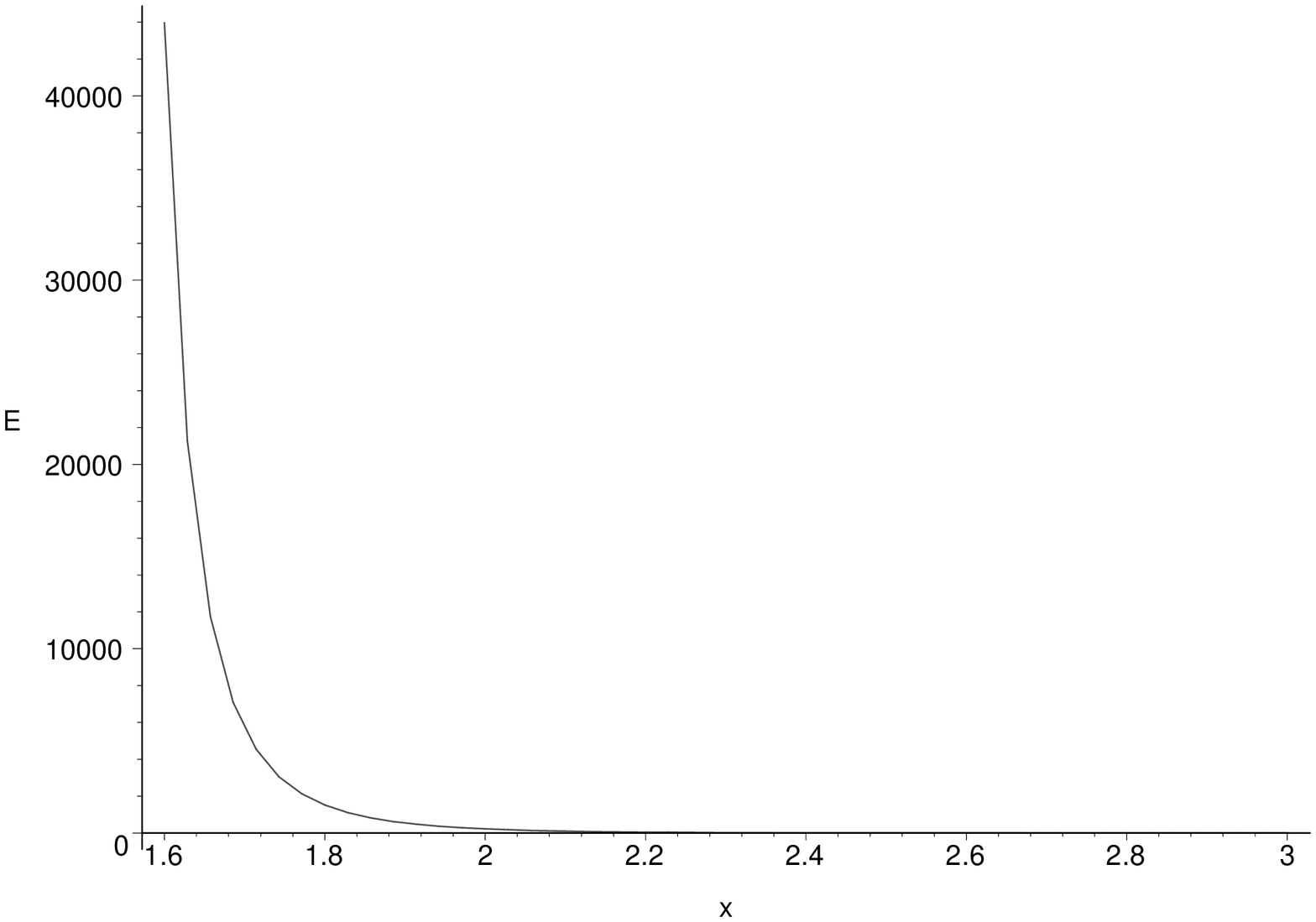,
height=4cm,width=5cm,angle=-90} \vspace{1cm} \caption{The parameter
$E(x,y)=E_0(x)$ corresponding to data considered in Example 2.}
\label{figure4}
\end{figure}

\textit{Example 2}. (Granville's egg curve) Consider $\partial
\Omega =\{(x,y):x^2y^2=(x-3)(1-x) ,\;x\in [1,3]\}$ (see Figure
\ref{figure3}) and the data given by
\[
w(x,y)=-y^2+{(x-3)(1-x)\over {x^2}}.
\]
The reduced ODE $2x(2x-3)E_0'(x)+2(x^4-4x+9)E_0(x)=x^4$ implies
\[
E(x,y)=E_0(x)={[C+E_1(x)]x^3\over {(2x-3)^{43/16}}} \exp
\left[{x\over {24}}(4x^2+9x+27)\right],
\]
which is defined for $x\neq 1.5$, where
$E_1'(x)=0.5(2x-3)^{27/16}\exp \left[ {x\over {24}}(4x^2+9x+27)
\right]$. Similarly, as in Example 1, the parameter cannot be
determined for $x=1.5$. See Figure \ref{figure4} for $E(1)=0.2$ and,
respectively, $E(3)=0.6$.

\vspace{0.5cm}

\begin{figure}[ht!]
\vspace{0.5cm} \epsfig{file=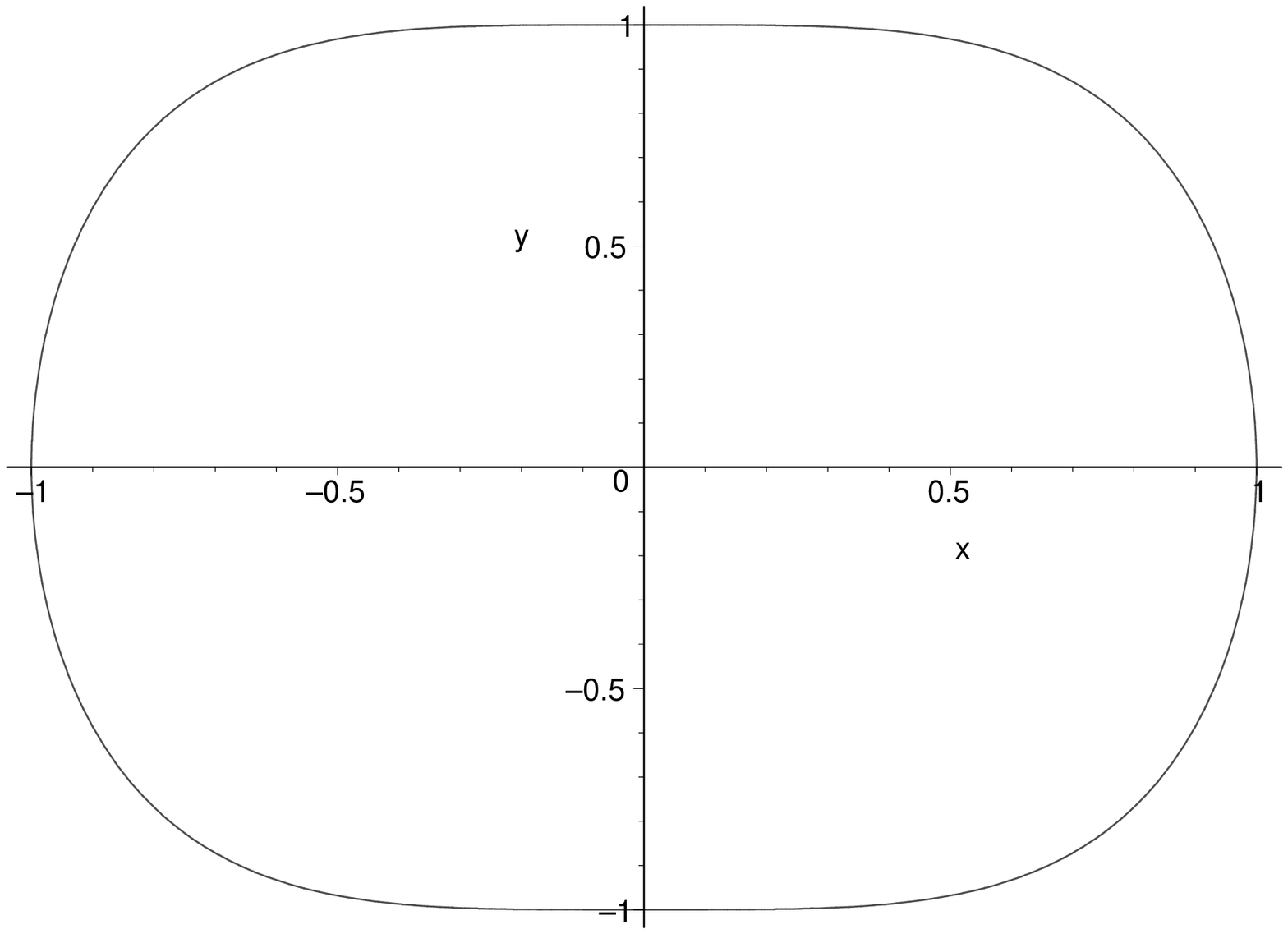,width=4cm, height=6cm,
angle=270} \hspace{0.75cm}
\epsfig{file=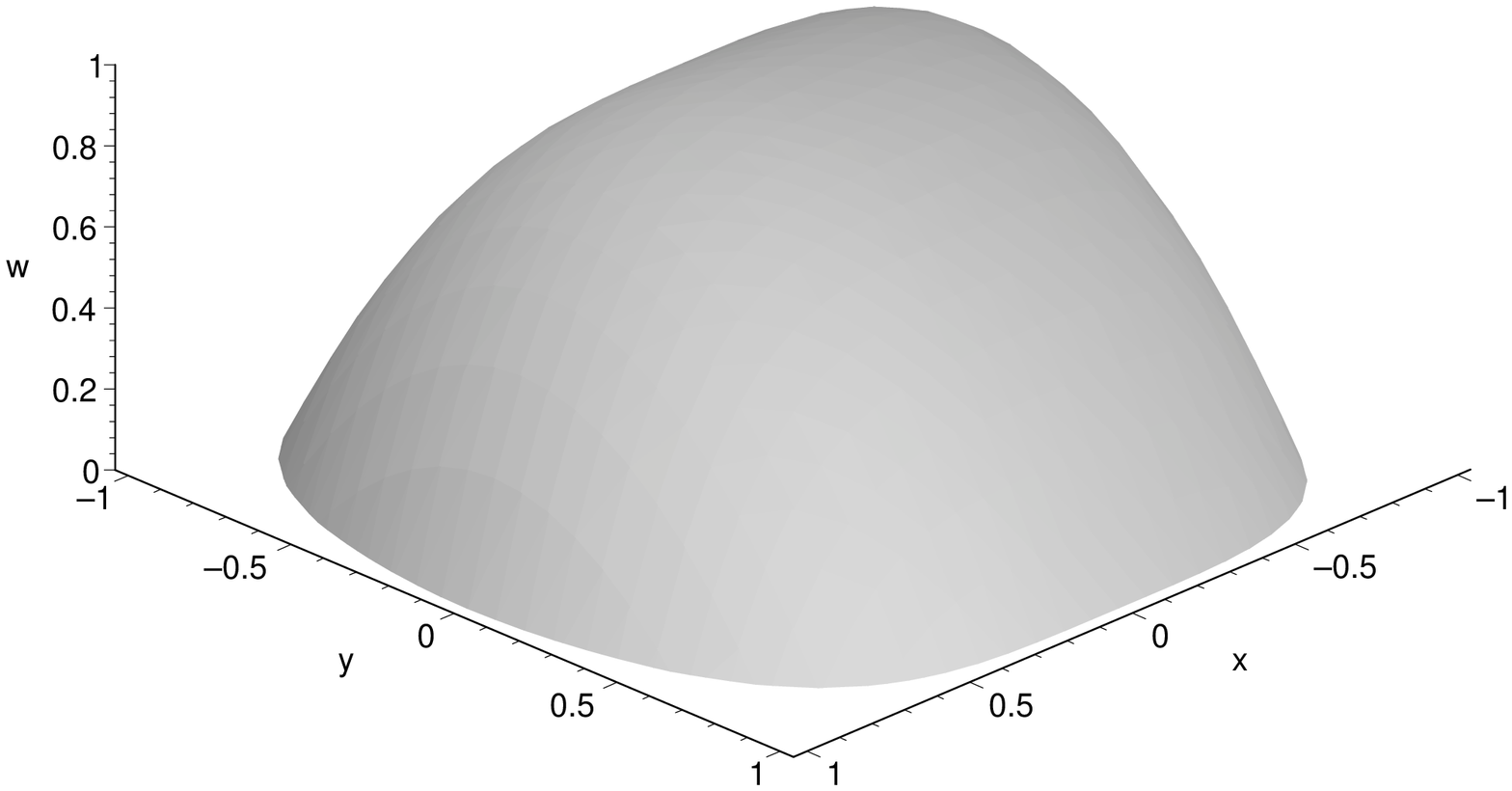,width=4cm,height=6cm,angle=-90} \vspace{1cm}
\caption{The boundary $\partial \Omega =\{(x,y):x^4+y^2=1 ,\;x\in
[-1,1]\}$ and the data $w(x,y)=1-x^4-y^2$ discussed in Example 3.}
\label{figure5}
\end{figure}

\begin{figure}[t!]
\vspace{0.5cm} \epsfig{file=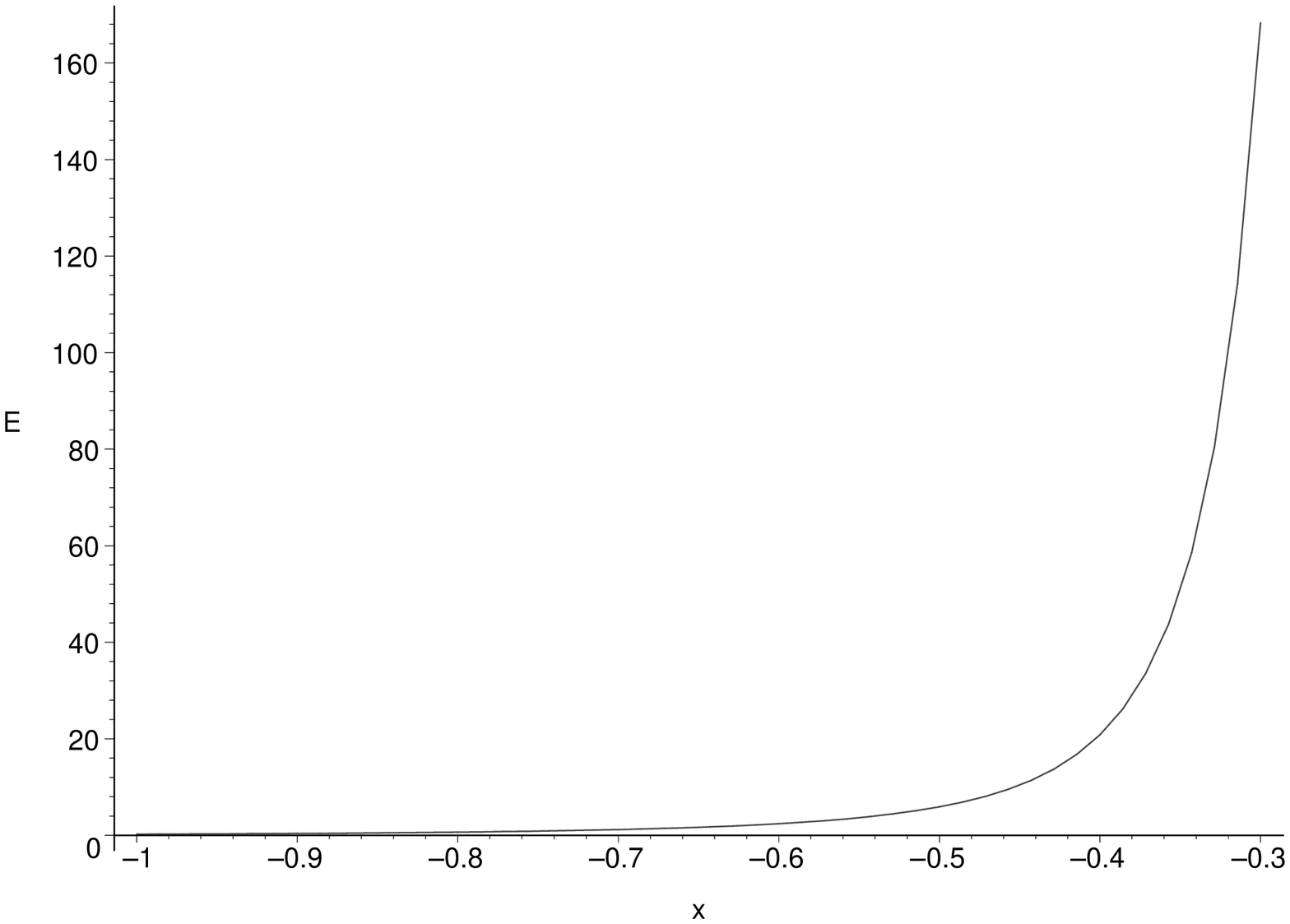,width=4cm, height=5cm,
angle=270} \hspace{1cm} \epsfig{file=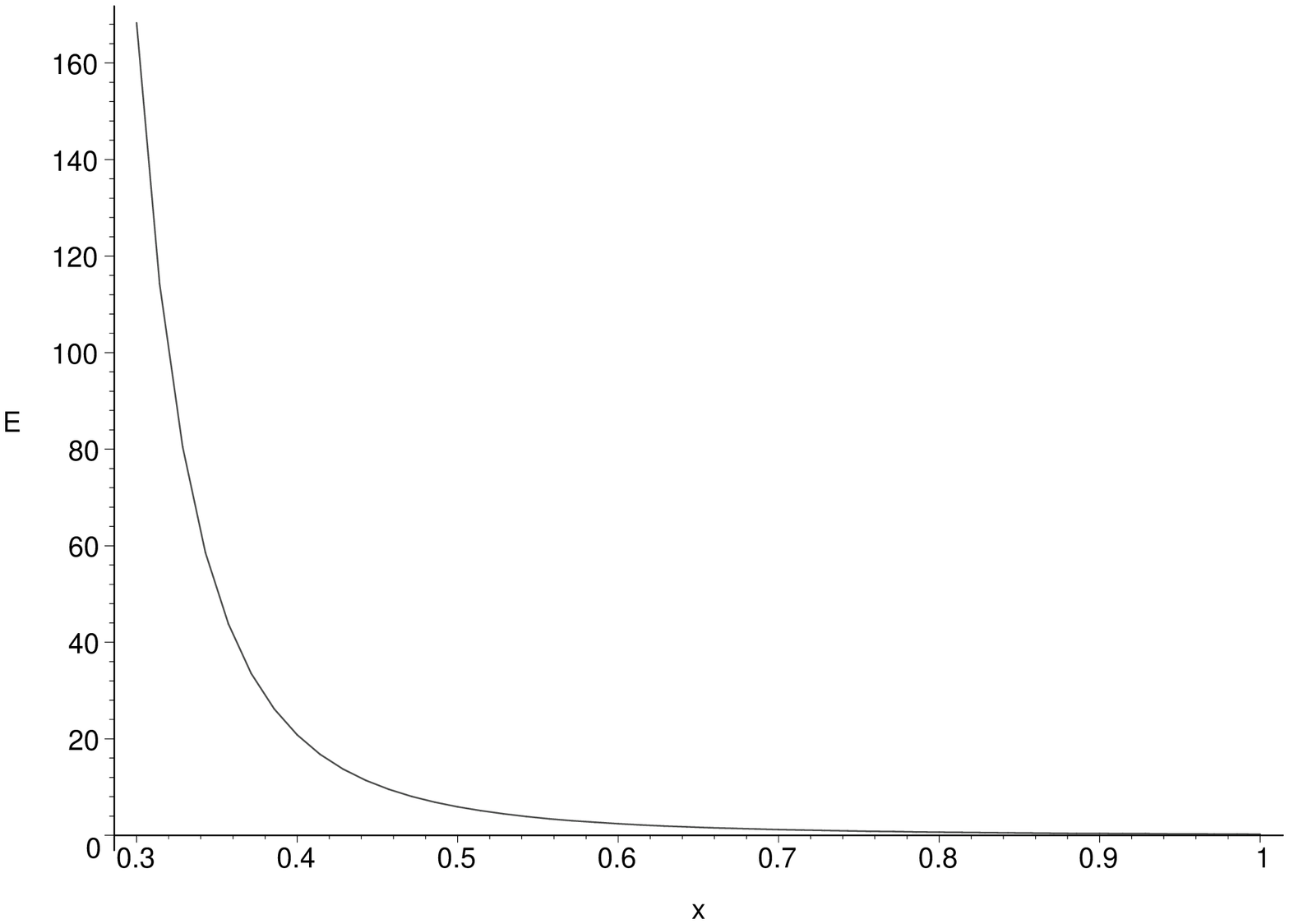,
height=5cm,width=4cm,angle=-90} \vspace{1cm} \caption{The parameter
$E(x,y)=E_0(x)$ corresponding to data discussed in Example 3.}
\label{figure6}
\end{figure}

\textit{Example 3}. (Generalized Lam\'{e} curve) For $\partial
\Omega =\{(x,y):x^4+y^2=1 ,\;x\in [-1,1]\}$ (see Figure
\ref{figure5}) and the following data
\[
w(x,y)=1-x^4-y^2,
\]
the PDE (\ref{eq2}) becomes $4x^3E_0'(x)+2(6x^2+1)=1$. It results
\[
E_0(x)=\frac{1}{4x^2}+\left[\frac{\sqrt{\pi}}{8x^3}
\hbox{erf}\left(\frac{1}{2x}\right)+\frac{C}{x^3}\right]\exp\left(
\frac{1}{4x^2} \right)
\]
that is defined for $x\neq 0$ (here $\hbox{erf}$ denote is the error
function). Therefore, the parameter cannot be estimated at $x=0$.
See Figure \ref{figure6} for the cases $E_0(-1)=0.25$, and,
respectively, $E_0(1)=0.25$.

\vspace{0.5cm}

\begin{figure}
\vspace{0.5cm} \epsfig{file=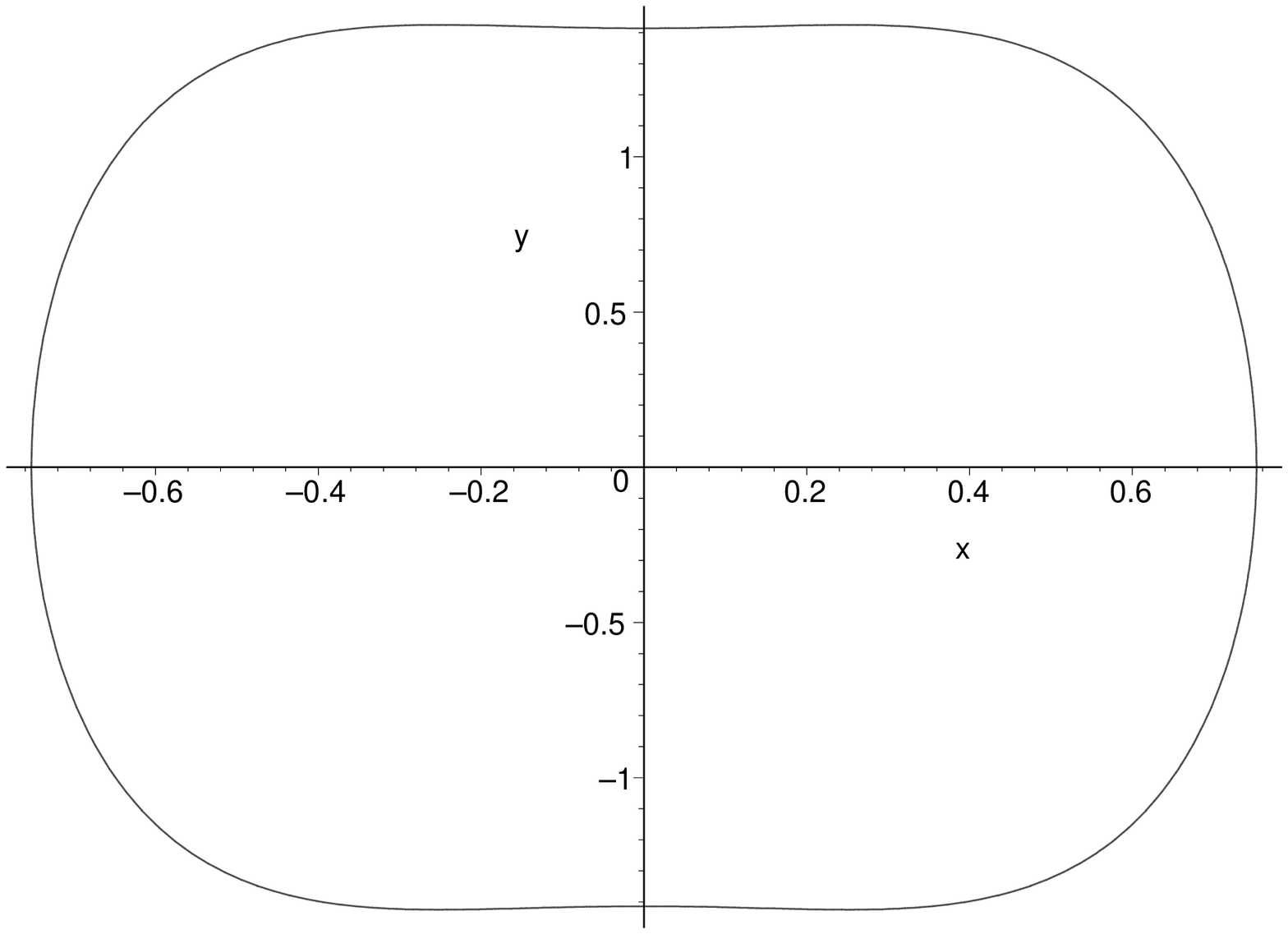,height=5cm,width=5cm,
angle=270} \hspace{1cm}
\epsfig{file=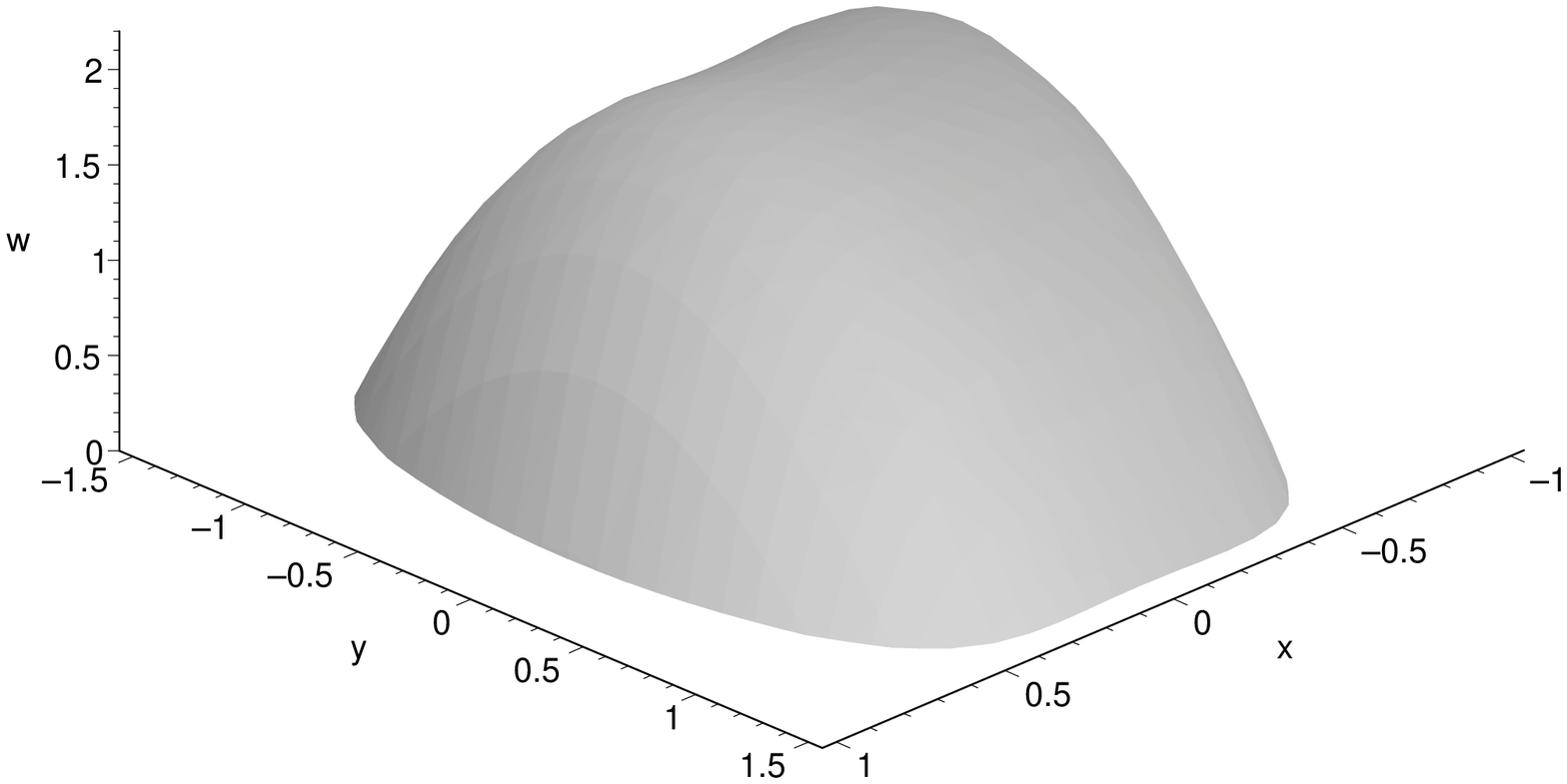,width=5cm,height=6cm,angle=-90} \vspace{1cm}
\caption{The boundary $\partial \Omega =\{(x,y):y^2 = -8x^4+ x^2+2,
x\in [-a,a]\}$, where $a=\sqrt{1+\sqrt{65}}/4$ and the data
$w(x,y)=2+x^2-8x^4-y^2$ considered in Example 4.}
\label{figure7}
\end{figure}

\begin{figure}
\vspace{0.5cm} \epsfig{file=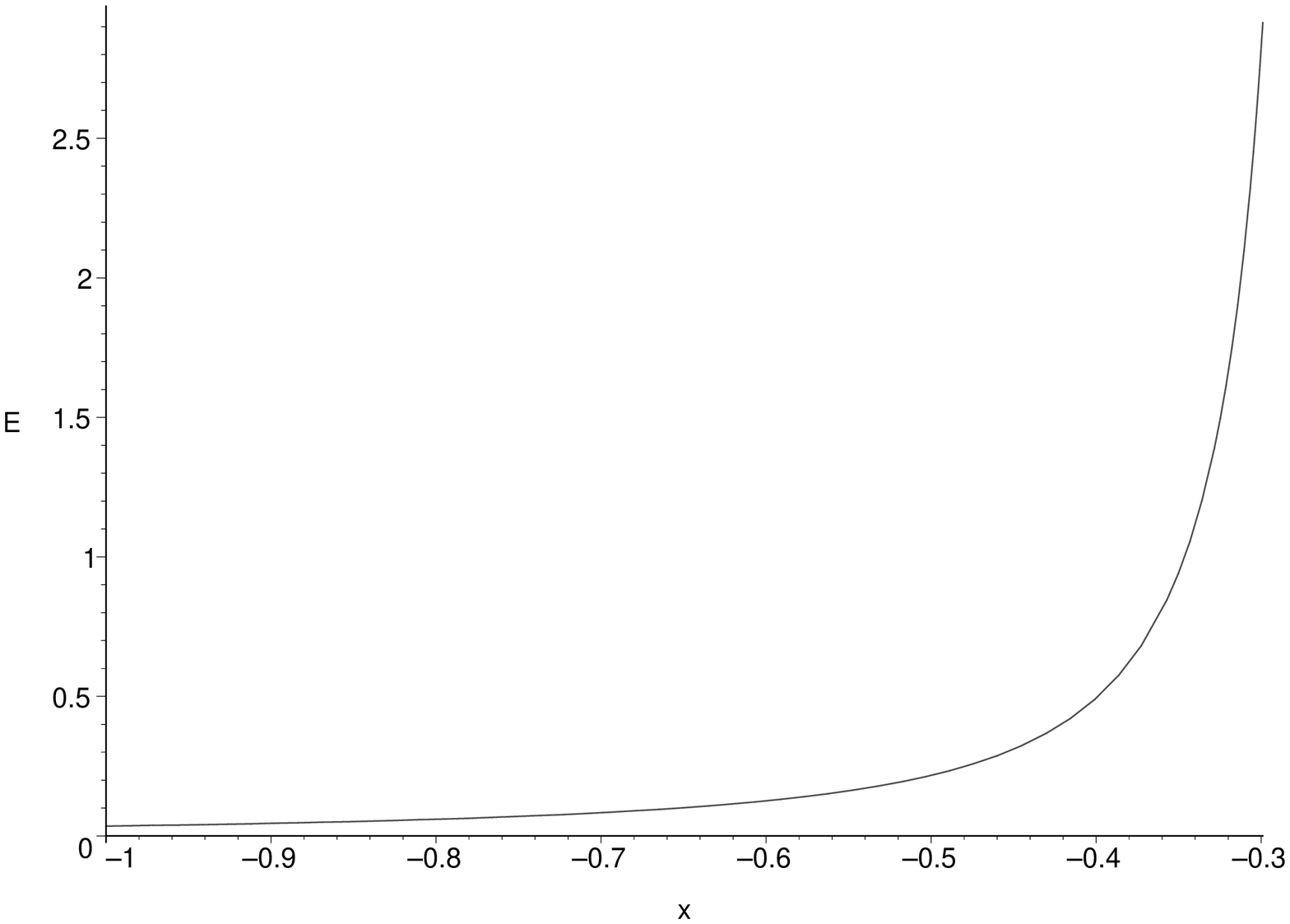,height=5cm,width=4cm,
angle=270} \epsfig{file=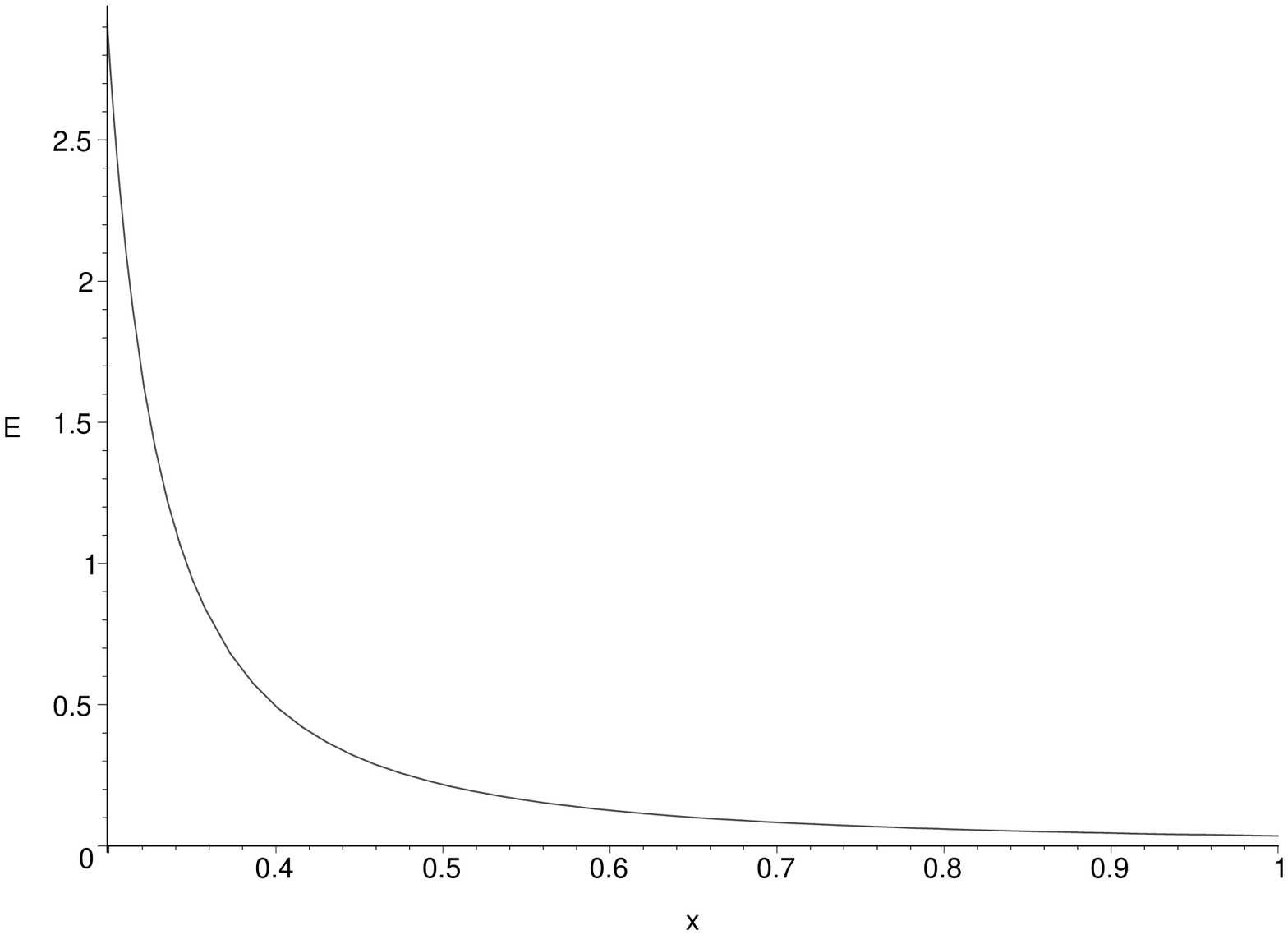,width=4cm,height=5cm,angle=-90}
\epsfig{file=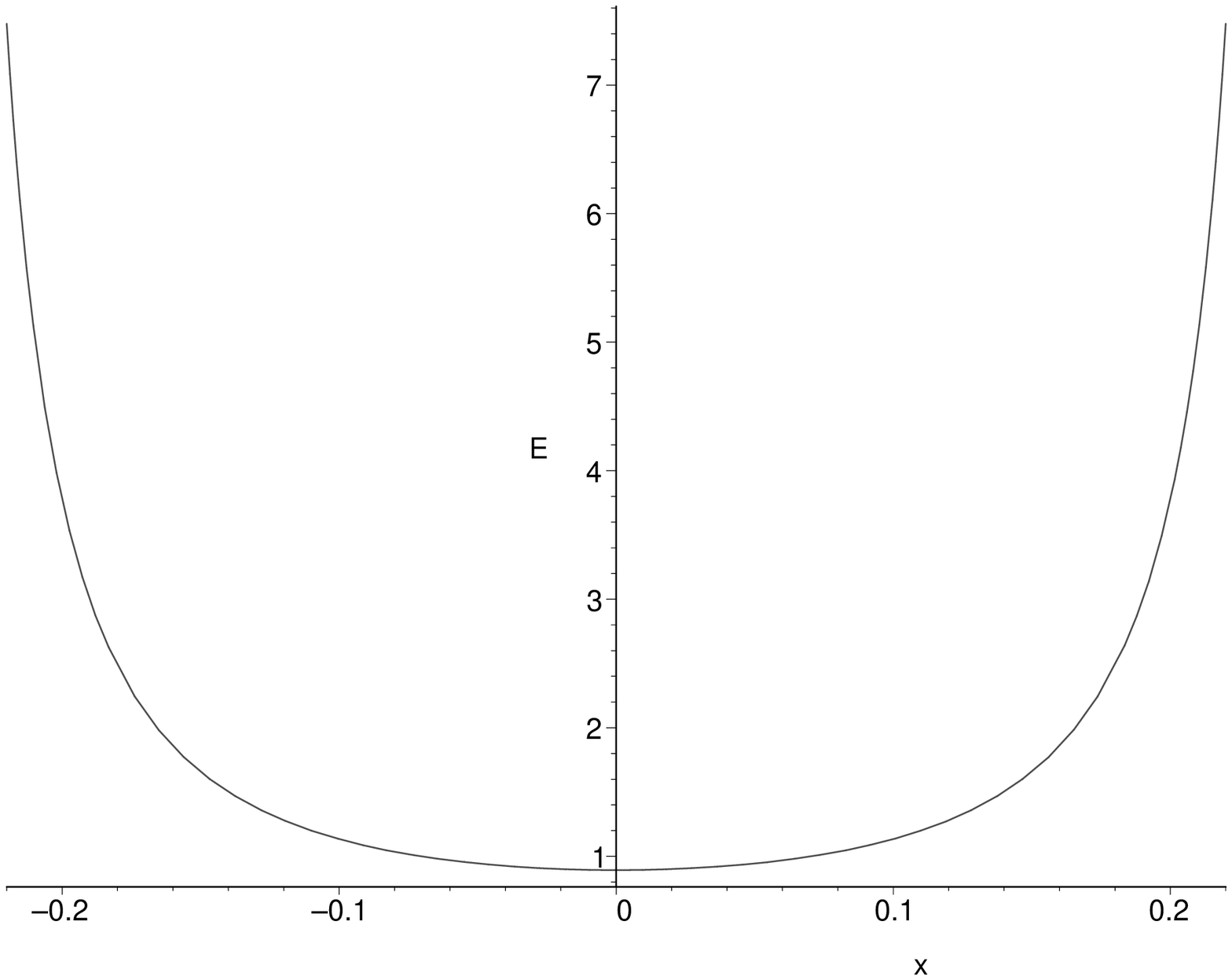,width=4cm,height=5cm,angle=-90} \vspace{1cm}
\caption{The parameter $E(x,y)=E_0(x)$ corresponding to the data
discussed in Example 4.}
\label{figure8}
\end{figure}

\textit{Example 4}. (Rounded corner rectangular domain) Suppose
$\partial \Omega =\{(x,y):y^2 = -8x^4+ x^2+2, x\in [-a,a]\}$, where
$a=\sqrt{1+\sqrt{65}}/4$ (see Figure \ref{figure7}). For the data
\[
w(x,y)=2+x^2-8x^4-y^2,
\]
we obtain the following reduced equation
$2x(16x^2-1)E_0'(x)+96x^2E_0(x)=1$. Thus,
\[
E_0(x)={{C_1+\sqrt{16x^2-1}+\arctan \left(\sqrt {16x^2-1}\right
)^{-1/2}}\over {2(16x^2-1)^{3/2}}}
\]
for $x\in [-a,-0.25)\cup (0.25,a]$, and
\[
E_0(x)={{C_2+\sqrt{1-16x^2}+\arctan \left(\sqrt {1-16x^2}\right
)^{-1/2}}\over {2(1-16x^2)^{3/2}}}
\]
for $x\in (-0.25,0.25)$. $E$ cannot be determined for $x=\pm 0.25$.
Figure \ref{figure8} shows the graph of the parameter for
$C_1=C_2=0$.

\vspace{0.5cm}

\section{Conclusion}
In this paper we  point out another systematic way of finding
classes of symmetry reductions related to parameter identification
problems of the form (\ref{ip}). Similar to the study in
\cite{bilacar}, we emphasize that the geometrical significance of
the nonlinearity occurring between the data and parameter in
(\ref{ip}) can be reflected by the group analysis tool. Seeking
different shapes for domains on which the dimension of the problem
can be reduced is not an easy task. Therefore, in this paper, we
discuss the nonclassical equivalence transformations related to
(\ref{eq2}). Briefly, to determine these transformations, the data
is considered as a dependent variable as well as the parameter, and
the nonclassical method (due to Bluman and Cole) is applied to the
studied equation. In general, for a known data $w$, one should check
the invariance of this function in (\ref{ISCEW2}), where $\xi$,
$\phi $, and $\psi $ are discussed in \S 3. Next, from the second
equation in (\ref{ISCEW2}), the form of the parameter should be
obtained in terms of the invariants of the symmetry reduction. At
the end, substituting $E$ and $w$ into (\ref{eq2}), the dimension of
the model should be reduced by one.

\textbf{Acknowledgements}. One of the authors, N. B\^{i}l\u{a} would
like to thank Prof. Dr. Heinz W. Engl at Johann Radon Institute for
Computational and Applied Mathematics (RICAM) who encouraged her to
apply the symmetry analysis to parameter identification problems.
She is also grateful to Dr. Philipp K\"{u}gler at RICAM and Dr.
Antonio Leit\~{a}o at Federal University of Santa Catarina, for many
stimulating and rewarding discussions during her stay at RICAM.


\begin{thebibliography}{99}


\bibitem{Ames} {\sc W. F. Ames}, {\em Nonlinear Partial
Differential Equations}, Academic Press, New York, 1967.

\bibitem{bilaniesen1} {\sc N. B\^{i}l\u{a} and J. Niesen}, {\em On a new procedure for finding nonclassical symmetries},
J. Symbolic Comput., 38, 6(2004), 1523-1533.

\bibitem{bilacar} {\sc N. B\^{i}l\u{a}}, {\em  Application of symmetry analysis to a PDE arising in the car windshield design},
SIAM Journal on Applied Mathematics, 65, 1(2004), pp.~113--130.


\bibitem{bluman1} {\sc G. W. Bluman and J. D. Cole}, {\em The general
similarity solutions of the heat equation}, J. Math. Mech., 18
(1969), pp.~1025--1042.

\bibitem{blumankumei} {\sc G. W. Bluman and S. Kumei}, {\em Symmetries and
Differential Equations}, Appl. Math. Sci.~81, Springer-Verlag, New
York, 1989.

\bibitem{BruGand} {\sc M. S. Bruz\'{o}na and M. L. Gandarias}, {\em Applying a new algorithm to derive nonclassical
symmetries}, Commun. in Nonlinear Sci. Numer. Simul., 13, 3(2008),
pp.~517--523.


\bibitem{desolv} {\sc J. Carminati and K. Vu}, {\em Symbolic
computation and differential equations: Lie symmetries}, J. Symbolic
Comput., 29 (2000), pp.~95--116.

\bibitem{Cheb1} {\sc E. S. Cheb-Terraba and A. D. Roched}, {\em Abel ODEs: Equivalence and
Integrable Classes}, Comput. Phys. Comm., 130, 1(2000),
pp.~204--231.


\bibitem{clarkson1} {\sc P. A. Clarkson and M. Kruskal}, {\em New
similarity reductions of the Boussinesq equation}, J. Math. Phys.,
30 (1989), pp.~2201--2213.


\bibitem{clarkson2} {\sc P. A. Clarkson and E. L. Mansfield},  {\em Algorithms for the nonclassical method of symmetry
reductions}, SIAM J. Appl. Math. 54(1994), pp.~1693--1719.

\bibitem{engl1} {\sc H. W. Engl, M. Hanke, and A. Neubauer}, {\em
Regularization of Inverse Problems}, Kluwer, Dordrecht, The
Netherlands, 1996.


\bibitem{Eriksson1} {\sc K. Eriksson, D. Estep, P. Hansbo, and C. Johnson}, {\em Computational Differential Equations},
Cambridge University Press, New York, 1996.

\bibitem{GGL1} {\sc G. Gambino, A. M. Greco, and M. C. Lombardo}, {\em A group analysis via weak
equivalence transformations for a model of tumor encapsulation}, J.
Phys A: Math. Gen. 37, 12 (2004), pp.~3835--3846.

\bibitem{hereman1} {\sc W. Hereman}, {\em Review of symbolic software for the
computation of Lie symmetries of differential equations}, Euromath
Bull., 1 (1994), pp.~45--82.

\bibitem{ibragimov} {\sc N. H. Ibragimov}, {\em CRC Handbook of Lie Group
Analysis of Differential Equations}, {\it Vol.} 1: {\em Symmetries,
Exact Solutions and Conservation Laws}, CRC Press, Boca Raton, 1994.


\bibitem{KL1} {\sc M. K\v{z}\'{i}\v{z}ek and L. Liu}, {\em Finite element approximation of a nonlinear
heat conduction problem in anisotropic media}, Comput. Methods Appl. Mech. Engrg., 157 (1998), p. 387—-397.


\bibitem{lie1} {\sc S. Lie}, {\em Gesammelte Abhandlungen}, Band
4, B.~G. Teubner, Leipzig, Germany, 1929, pp.~320--384.

\bibitem{meleshko1} {\sc S. V. Meleshko}, {\em Generalization of the equivalence transformations}, J. Nonlinear
Math. Phys., 3, 1-2(1996), pp.~170--174.

\bibitem{olverbook} {\sc P. J. Olver}, {\em Applications of Lie Groups to Differential
Equations}, Grad. Texts in Math.~107, Springer-Verlag, New York,
1986.

\bibitem{ovsyannikov1} {\sc L. V. Ovsyannikov}, {\em Group Analysis of
Differential Equations}, W.~F. Ames, trans., Academic Press, New
York, 1982.

\bibitem{handbookPZ} {\sc A. D. Polyanin and V. F. Zaitsev}, {\em Handbook of Exact Solutions for
Ordinary Differential Equations}, Chapman \& Hall/CRC Press, Boca
Raton, Second Edition, 2003.


\bibitem{handbookPZ2} {\sc A. D. Polyanin and V. F. Zaitsev}, {\em Handbook of Nonlinear Partial Differential Equations},
Chapman \& Hall/CRC Press, Boca Raton, 2004.

\bibitem{torrisi1} {\sc V. Romano and M. Torrisi}, {\em Application of weak equivalence
transformations to a group analysis of a drift-diffusion model}, J.
Phys. A: Math. Gen. 32 (1999), pp.~7953--7963.

\bibitem{torrisi3} {\sc M. Torrisi and R. Tracin\`{a}}, {\em Equivalence transformations and symmetries
for a heat conduction model}, Int. J. of Non-Linear Mechanics, 33
(1998), pp.~473--487.

\bibitem{zwillinger} {\sc D. Zwillinger}, {\em Handbook of Differential Equations}, Third Edition,
Academic Press, Boston, 1997.

\end{thebibliography}
\end{document}